\documentclass[sn-mathphys]{sn-jnl}

\jyear{2021}%

\theoremstyle{thmstyleone}%
\newtheorem{theorem}{Theorem}
\newtheorem{proposition}[theorem]{Proposition}%

\theoremstyle{thmstyletwo}%
\newtheorem{example}{Example}%
\newtheorem{remark}{Remark}%

\theoremstyle{thmstylethree}%
\newtheorem{definition}{Definition}%
\usepackage{wasysym}
\usepackage{empheq}
\usepackage{url}

\usepackage{float,multirow,booktabs}\usepackage{array}
\usepackage[novbox]{pdfsync} \usepackage{relsize,exscale}
\usepackage[normalem]{ulem}
\usepackage{hyperref}
\usepackage{cleveref}
\usepackage{upgreek}
\usepackage{color}
\usepackage{xcolor}
\usepackage{subcaption}

\usepackage{graphicx}
\usepackage{epsfig}
\graphicspath{{figures/}}

\usepackage{geometry}

\usepackage{vmargin}
\setmarginsrb{.8in}{.6in}{.8in}{.4in}{0pt}{.3in}{0pt}{.3in}

\newcommand*{\QEDB}{\hfill\ensuremath{\square}}

\newcommand{\RN}[1]{\textup{\uppercase\expandafter{\romannumeral#1}}}

\def\L{\Lambda}
\newcommand{\mV}{\mathcal{V}}

\newcommand{\N}{\mathbb{N}}
 \def\satd{satisfied}  \def\itc{it is easy to check that }\def\wmw{we may write}

\newcommand{\n}{\mathsf n}

\def\How{However, }\def\PF{Perron-Frobenius}
\def\no{\nonumber}\def\bep{\begin{pmatrix}} \def\eep{\end{pmatrix}}
 \def\DFE{disease free equilibrium}\def\bn{\bff \nu}

\def\l{{\lambda}}

\newcommand{\realnneg}{\mathbb{R}_{\geq 0}}

\providecommand{\pp}[1]{\left[#1\right]} 
\providecommand{\pr}[1]{\left(#1\right)} 

{Lemma}
\newtheorem{corollary}
{Corollary}
\def\beC{\begin{corollary}
  }\def\eeC{\end{corollary}}
{Conjecture}

\newcommand*{\QED}{\hfill\ensuremath{\blacksquare}}

\renewcommand{\theta}{\vartheta}

\renewcommand{\thefootnote}{\fnsymbol{footnote}}

\numberwithin{equation}{section}

\def\bc{\begin{cases}
  }      \def\mF{\mathcal F}
\def\ec{\end{cases}}  
  \def\qu{\quad} \def\for{\forall}

  \newcommand{\beq}{\begin{eqnarray}
    }
\def\eeq{\end{eqnarray}}
   \newcommand{\be}[1]{\begin{equation}\label{#1}}
\newcommand{\ee}{\end{equation}}
\def\bea{\begin{eqnarray*}}\def\ssec{\subsection}
  \def\Prf{{\bf Proof:} }  
\def\eea{\end{eqnarray*}} \def\elt{element} \def\la{\label}\def\fe{for example }    
   \def\Mp{More precisely, } \def\sats{satisfies}  \def\saty{satisfy}       

\def\I{\infty} \def\Eq{\Leftrightarrow}
  \def\T{\widetilde}
\def\PH{phase-type }  
\def\BEN{\begin{enumerate}}  \def\BI{\begin{itemize}}
\def\EEN{\end{enumerate}}   \def\EI{\end{itemize}} \def\im{\item} \def\Lra{\Longrightarrow}  \def\eqr{\eqref}  
\def\no{\nonumber} 
 \def\nR{\mathcal R}
\def\mR{\mathcal R_0} \def\mI{\mathcal I}
\def\g{\gamma}     \def\b{\beta}
   
 \def\eig{eigenvalue}\def\eigv{eigenvector}
   
\def\Fr{Furthermore, }
    
   \def\Mp{More precisely, } 
   \def\wrt{with respect to }
  \def\mbw{may be written as }\def\resp{respectively}   
 \def\eqr{\eqref}  
\def\wk{well-known}

\def\l{\; \mathsf l}
\def\fr{\frac} \def\im{\item}
\newcommand{\s}{\;\mathsf s}
\renewcommand{\i}{\;\mathsf i}
\def\m0{{\mathcal R}_0}

\def\eeD{\end{definition}} \def\beD{\begin{definition}}
\def\beR{\begin{remark}} \def\eeR{\end{remark}}
\def\mD{\mathcal D }
       \def\bb{\bff b}  
 \def\vi{\vec \i}  \def\vz{\vec 0}\def\vg{\vec \g} \def\bg{\bff \gamma} \def\sd{s_{dfe}}
    \def\va{\vec \alpha } 
     \def\vr{\; \vec {\mathsf r}}\newcommand{\bff}[1]{{\mbox{\boldmath$#1$}}}
    \def\beP{\begin{proposition}} \def\eeP{\end{proposition}}
      
    \def\it{immunity threshold}
    
    \long\def\symbolfootnote[#1]#2{
\begingroup
\def\thefootnote{\fnsymbol{footnote}}\footnote[#1]{#2}
\endgroup}
\def\fn{\symbolfootnote}
\def\corr{corresponding }
\def\l{\lambda}
\def\com{compartment}\def\nne{nonnegative}
\newcommand{\rt}{\right}

\def\D{\Delta}

\def\a{\; \mathsf a}   
\def\brn{basic reproduction number}\def\ngm{next generation matrix}
\newcommand{\red}{\textcolor[rgb]{1.00,0.00,0.00}}

\renewcommand{\e}{\;\mathsf e}
\renewcommand{\r}{\;\mathsf r}
 
\begin{document}

\title[On varying population matrix-SIR Arino models]{On  matrix-SIR Arino models with linear birth rate, loss of immunity,  disease and vaccination fatalities, and their approximations}


\author*[1]{\fnm{Florin} \sur{Avram}}\email{florin.avram@univ-pau.fr}

\author[2]{\fnm{Rim} \sur{Adenane}}\email{rim.adenane@uit.ac.ma}
\equalcont{These authors contributed equally to this work.}

\author[3]{\fnm{Lasko} \sur{Basnarkov}}\email{lasko.basnarkov@finki.ukim.mk}
\equalcont{These authors contributed equally to this work.}

\author*[4]{\fnm{Gianluca} \sur{Bianchin}}\email{gianluca.bianchin@colorado.edu}
\equalcont{These authors contributed equally to this work.}

\author[5]{\fnm{Dan } \sur{Goreac}}\email{dan.goreac@u-pem.fr}
\equalcont{These authors contributed equally to this work.}

\author[6]{\fnm{Andrei  } \sur{Halanay}}\email{andrei.halanay@upb.ro}
\equalcont{These authors contributed equally to this work.}

\affil*[1]{\orgdiv{Laboratoire de Math\'ematiques Appliqu\'es}, \orgname{Universit\'e de Pau}, \orgaddress{ \city{ Pau}, \postcode{64000},  \country{France}}}

\affil[2]{\orgdiv{Department of Mathematics}, \orgname{Ibn-Tofail University}, \orgaddress{ \city{Kenitra}, \postcode{14000},  \country{Morocco}}}

\affil[3]{\orgdiv{Faculty of Computer Science and Engineering}, \orgname{Ss. Cyril and Methodius University in Skopje}, \orgaddress{ \city{Skopje},  \country{North Macedonia}}}

\affil*[4]{\orgdiv{Department of Electrical, Computer, and Energy Engineering}, \orgname{University of Colorado Boulder}, \orgaddress{ \city{Boulder}, \postcode{CO 80309},  \country{United States}}}

\affil[5]{\orgdiv{School of Mathematics and Statistics}, \orgname{ Shandong University}, \orgaddress{\postcode{ 264209}, \city{Weihai},  \country{China}}}

\affil[5]{\orgdiv{LAMA, Univ Gustave Eiffel}, \orgname{ UPEM, Univ Paris Est Creteil}, \orgaddress{\postcode{ CNRS, F-77447}, \city{Marne-la-Valleée},  \country{France}}}

\affil[6]{\orgdiv{Department of mathematics and informatics}, \orgname{ Polytechnic University of Bucharest}, \orgaddress{ \city{Bucharest},  \country{Romania}}}


\abstract{In this work we study the stability properties of the
equilibrium points of deterministic epidemic models with nonconstant
population size.
Models with nonconstant population have been studied in the past only
in particular cases, two of which we review and combine.  Our main
result shows that for simple ``matrix epidemic models" introduced
in~\cite{Arino}, an explicit general formula for the reproduction
number $\nR$ and the corresponding ``weak  stability alternative"
\cite[Thm 1]{Van08} still holds, under small modifications, for models
with nonconstant population size, and even when the model allows for
vaccination and loss of immunity.
The importance of this result is clear once we note that  the  models
of \cite{Arino} include a large number of  viral and bacterial models
of epidemic propagation, including for example the totality of
homogeneous COVID-19 models.
To better understand the nature of the result, we emphasize that the
models proposed in \cite{Arino} and considered here are  extensions of
the SIR-PH model \cite{Riano},  which is  essentially characterized by
a \PH\ distribution $(\va, A)$ that models transitions between the
``disease/infectious compartments".
In these cases, the reproduction number  $\nR$ and a certain Lyapunov
function for the \DFE\  are explicitly expressible in terms of
$(\va, A)$.
Not surprisingly, accounting for varying demography, loss of immunity,
and  vaccinations lead to  several challenges.
One of the most important is that a varying population size  leads to
multiple endemic equilibrium points: this is in contrast with
``classic models,'' which in general admit unique disease-free and
endemic equilibria.
As a special case of our analysis, we consider a ``first approximation"
(FA) of our model, which coincides with the constant-demography model
often studied in the literature, and for which more explicit results
are available. \Fr we propose a second heuristic approximation named
``intermediate approximation" (IA).
We hope that more light  on  varying population models with loss of
immunity  and  vaccination, which have been largely avoided until now
-- see though \cite{busenberg1990analysis,busenberg1993method,Der,Green97,LiGraef,SunHsieh,Raz,Arinovar}
-- will emerge in the future.
}




\maketitle

\newpage

\section{Introduction}\label{sec1}
In this paper, we extend -- to epidemic models with
varying population size, vaccinations, and loss of immunity -- an
outstanding  formula derived in \cite{Arino}, which expresses the
reproduction number of a large class of deterministic  epidemic models
as a quadratic form \wrt\ certain input parameters \eqr{nRV}. We
illustrate why this formula is remarkable through a brief  historical overview.

\smallskip
{\bf Deterministic mathematical epidemiology.} {Deterministic
mathematical models~\cite{Mart,Thieme,Chavez} have been widely
adopted by epidemiologists to model the spread of diseases worldwide,
including  the  Bombay plague in 1905--06 \cite{KeMcK},   measles,
smallpox, chickenpox, mumps, typhoid fever and diphtheria (see \fe~\cite{earn2008light}), and~recently COVID-19 ( see \fe\  \cite{Schaback,bacaer2020modele,Ketch,Charp,Djidjou,Sofonea,alvarez2020simple,
horstmeyer2020balancing,
di2020optimal,Franco,baker2020reactive,caulkins2020long,caulkins2021optimal}, to~cite just a few representatives of a huge literature).}
Analytic work in this area appears to have reached the ''dimensionality
barrier," thus remaining limited to models with up to three epidemic
states~\cite{busenberg1993method,Der}.

\smallskip
{\bf Stochastic  and deterministic  models}. Stochastic models are
the most-natural way to approach epidemic modeling since they
inherently capture stochastic mixing within
populations~ \cite{bartlett2020deterministic,baxendale2011sustained}.
We remark that many deterministic models can be obtained as limiting
cases of stochastic models~\cite{kurtz1978strong,Brit}.
Notice, however, that although multiple stochastic models may be
adequate to model a certain epidemic outbreak, these may yield
different deterministic limits -- see \fe\ \cite{naasell2013influence},
who discusses four stochastic models for a particular study case.
The correspondence between  stochastic and  deterministic models being
a delicate point, it  implies that well-defined  concepts for the first
class are more challenging to capture or interpret for the second. To
further explain  this point and our motivation, next we give an
eye-bird's view of mathematical epidemiology.

The deterministic epidemics literature may be divided roughly into three streams.
\BEN \im Models with constant total population size.
These models are in general easier to study, often admit a single
endemic equilibrium point,  and typically obey the ``$\mR$ alternative" -- see the next class of models.
 \How\ since death is an essential factor of epidemics, the assumption of constant population size (clearly an approximation that holds in the short term or for very-large populations)  deserves some  comment%
\fn[1]{One  rigorous justification for deterministic constant population  epidemiological models comes from slow-fast analysis \cite{Kuehn}.   This is best understood for models with demography (birth,  death), which happen typically on a  slower scale than the infectious  phenomena. Here there is a natural partition of the compartments
into a vector $\vi(t) \in \mathbb{R}_+^-$ of disease/infectious \com s (asymptomatic, infectious hospitalized,etc).  These interact (quickly) with the other input classes like the susceptibles and output classes like the recovered and dead.}. 

\im Models with constant birth rate. These models include the previous
class, and preserve some of its features, such as uniqueness of the
endemic fixed point. In general cases, their stability properties can
be studied via the  \ngm\ approach -- see \cite{de2019some} for a
recent  nice review of several stability results for  this class  of
models. Finally, we notice that these models correspond to limiting cases of stochastic models with emigration.

\im  Finally, we arrive to the object of our paper:  models with   linear birth-rate (corresponding to constant birth-rate per capita in
the analog stochastic model), varying total population, and
proportionate mixing. This stream of literature  precedes the \ngm\ revolution \cite{busenberg1990analysis,busenberg1993method}, and  reveals the possibility of bi-stability when $\mR <1$ (absent from the previous models), even in simple examples\fn[5]{Thus, for an initial
number of infectives high, the trajectory may  lie in the
basin of attraction  of a stable
endemic  point instead of being eradicated. The discrepancy with what is expected from the corresponding stochastic model suggests that in this range the deterministic model is inappropriate}  \cite{Der}.
Despite further remarkable works  on particular cases -- see, for example, \cite{Green97,LiGraef,Raz,li2002qualitative,SunHsieh,Arinovar
} ({which proposed a direct stability analyses as opposed to the \ngm\
approach used here}), the  literature on models with varying total population, unlike the  two preceding streams, has not yet reached  general  results.

\EEN

Despite their importance in describing epidemics evolving over long
periods of time or affecting small-size populations, the previous
discussion suggests that epidemics model with varying population size
and linear birth-rate have not not received sufficient attention --
even in the simple case of Arino models with linear forces of infection (or bilinear incidence) \cite{Arino,feng2000endemic,Feng,Riano,AAK} --
thus motivating our study.
As we will show in this work, these are models to which the \ngm\
approach can be applied, and for which the \ngm\ has rank $1$,
resulting  in an explicit formula for the basic reproduction number
$\mR$, provided in \cite{Arino} for the case of constant population
size\fn[5]{Unfortunately, this formula is not known enough, and particular cases of it  are being  reproved in numerous recent papers.
In fact,  sometimes several papers reprove the same particular case, due to  the confusion caused  by the lack of a common notation style;  we have a proposition below to remedy that.}.

We remark that Arino \cite{Arino} calls the class of models
considered in this work ``simple models" \cite{arino2020simple}. These
models were introduced in Arino's joint work \cite{Arino} with several
of the founders of mathematical epidemiology: Brauer,  van den Driessche,  Watmough, and Wu.
We took the liberty of renaming this class, since we believe that
``simple models'' is too generic, while ``Arino-Brauer-van den Driessche-Watmough-Wu'' is too long.  We hope to  show here and in further work that for this class $\mR$ and other important features
continue  to have  simple formulas  in terms of the intervening  matrices, and of the vaccination and loss of immunity parameters\fn[4]{Note also that quadratic models have been found useful in \cite{buonomo2015note} for determining the direction of transcritical bifurcations.}.

\smallskip
{\bf Related works}. We next summarize some works that are related to
the model considered in this work.
\BEN \im
The simplest class of ``matrix-Arino" epidemics are the ``SIR-PH
models" \cite{Riano,AAK}, in which the infected class is replaced by a
vector of disease classes (here denoted by $I$), there are several
output  classes (here denoted by $R$),  and there is only one input
class (here denoted by $S$).    Note that this particularly  simple model has the  probabilistic interpretation  of a   SIR model where the exponential infection time has been replaced by a PH-type $(\va, A)$ distribution \cite{Riano}.\fn[5]{One may similarly replace the exponential latency  time in class E of SEIR by a PH-type \cite{Hurtado}, and similarly for all the infectious classes, but this is finally
unnecessary, since all the infectious classes may be grouped together in one group, whose \PH will be determined by those of the components (via Kronecker product operations). It suffices therefore to let $(\va, A)$ denote the \PH of all the interactions
between the disease classes.}
Note that many of the models recently adopted to model the COVID-19
outbreak belong  to this class.

\im
Unfortunately, the SIR-PH model precludes important interactions
between $S$ and $R$ like loss of immunity and  vaccination.  We
introduce therefore here a  class of models, named SIR/V+S, which allows
individuals to transfer from the $S$ class to the $R$ class (thus
accounting for vaccinations) and from the $R$ class to the $S$ class
(thus accounting for loss of immunity).
   \EEN

{\bf Contributions}. The contribution of this work is twofold. First,
we propose a compartmental model for epidemic transmission that
accounts for varying population size, vaccinations, and loss of
immunity. This model includes three well-studied models as subcases:
(i) general varying population epidemic models, (ii) the well-studied
``first order" approximations (FA), often adopted in the
literature (which is recovered by ignoring certain quadratic terms
from the more-general model proposed here), and (iii) a new model,
named  ``intermediate approximations" (IA), introduced here for the
first time  (obtained by neglect the terms which are quadratic \wrt\
the disease/infectious compartments).

Second, we characterize the equilibrium points of the proposed model
and we study their stability properties. Our analysis builds upon and
extends the \ngm\ approach proposed in \cite{de2019some} for models
with constant population size.
Our results show that the formula for the reproduction number proposed
in \cite{Arino} extends to models with varying population (with small
modifications) and that the reproduction number still characterizes the
stability of the equilibrium points.

{\bf  Organization}. This paper is organized as follows. In Section 2
we present some preliminary background on epidemic models: the basic reproduction number
and the \ngm\ method.
In Section \ref{s:Feng} we introduce matrix Arino models with demography, vaccination and loss of immunity. Section \ref{s:S} gives stability results for the SIR-PH model
-- see Propositions \ref{p:brn}, \ref{p:ee}, \ref{p:IA}. Section 5 concludes the paper.

\section{Preliminaries\la{s:pr}}
\ssec{What is a deterministic epidemic model?}
To put in perspective our work, we would like to start by a  definition of deterministic epidemic  models, lifted from \cite{KamSal}.

\beD
A deterministic epidemic  model is a  dynamical system with two types of variables $\vec x(t):=(\vi(t),\vec z(t))\in {\mathbb{R}}_{+}^{N}$, where
\BEN \im
$\vi(t)$ model the number (or  density) in
 different compartments of  infected individuals (i.e. latent,
infectious, hospitalized, etc) which should ideally disappear eventually if the epidemic ever ends;

\im $\vec z(t)$ model  numbers (or  densities) in
  compartments of individuals
who are not infected (i.e. susceptibles, immunes,
recovered individuals, etc).
\EEN

 The system must admit  an  equilibrium called \DFE\ (DFE), and hence a ``quasi-triangular" linearization the form
\beq \la{dyn} \vi'(t) &&=\vi(t) A_{i,i}(\vec x(t)), \vec x(t)\in \mD \subset{\mathbb{R}}_{+}^{N}\\
\vec z'(t)&&=\vi(t) A_{z,i}(\vec x(t)) + (\vec z(t)-\vec z_{dfe})  A_{z,z}(\vec x(t)),\no\eeq
where $\mD$  is some forward-invariant subset,
where ``quasi-triangular" refers to the fact that the functions $A_{i,i}, A_{z,i},A_{z,z}$ depend on
all the variables $\vec x(t)$, and  where $N$ is the dimension of $\vec x(t)$.
\QEDB
\eeD

As shown in  \cite{KamSal}, any epidemic model admitting an equilibrium
point  $(\vz, \vec z_{dfe})$ admits the representation \eqr{dyn}, under
suitable smoothness assumptions.
In what follows, we will call the point $x_{dfe}=(\vz, \vec z_{dfe})$ a
\DFE\  (DFE).

\beR
Note that the  essential feature of \eqr{dyn} is the ``factorization of the disease equations". 
\QEDB
\eeR

\ssec{The basic reproduction number  $\mR$}
{The basic reproduction number  or ``net reproduction rate" $\m0$  is a pillar concept in demography, population dynamics, branching processes and mathematical epidemiology -- see the introduction of the book \cite{bacaer2021mathematiques}.}
One of the central objectives of these fields
is to study the stability of DFE, i.e. the conditions which ensure eradicating the sickness (or a part of the population in population dynamics). It was discovered in  simple models, that this amounts  to verifying that a famous threshold parameter called basic reproduction number   is less than $1$.

 \BEN \im The notation $\m0$ was  first introduced by the father
of mathematical demography Lotka \cite{lotka1939analyse,dietz1993estimation}. In epidemiology, the basic reproduction number models
 the expected number of secondary cases which one infected case would produce in a  homogeneous,
completely susceptible stochastic population, in the next generation.  As well known in the simplest setup of branching  process, having this  parameter smaller than $1$ makes extinction sure.  The relation to epidemiology is   that an epidemics is well approximated by a branching process at its inception, a fact which goes back to  Bartlett and Kendall.
\im With more infectious classes, one deals at inception with  multi-class branching processes, and stability holds when the \PF\ eigenvalue of the   ``{\bf \ngm\ }" (NGM)  of means  is smaller than $1$.

  \im   For  deterministic epidemic models, it seems at first that the basic reproduction number  $\m0$ is lost, since the  generations disappear in this setup -- but see  \cite[Ch. 3]{bacaer2021mathematiques}, who recalls a method to introduce generations which goes back to Lotka, and which is reminiscent of the iterative Lotka-Volterra approach of solving integro-differential equations. At the end of the tunnel, a unified method for defining $\m0$ emerged only much later, via
      the ``next generation matrix"  approach \cite{diekmann1990definition,Van,Van08,Diek,perasso2018introduction}.
The final result is that local stability of the \DFE\ holds iff the spectral radius of a certain matrix called ``\ngm",  which depends only on a set  of  ``infectious  compartments" $\vi$ (which we aim to   reduce to $0$), is less than one.
This  approach works provided that certain  assumptions listed below hold\fn[4]{And so $\m0$ is undefined when these assumptions are not \satd.}.
\EEN

      \BEN \item[(C1)] The foremost assumption is that the  disease-free  equilibrium $(\vec 0, \vec z_{dfe})$ is
{\em unique and locally asymptotically stable within the disease-free space} $\vi=0$, meaning that all solutions of   $$\vec z'(t)=(\vec z(t)-\vec z_{dfe})  A_{z,z}(\vz,\vec z(t)), \qu \vec z(0)=\vec z_0$$
 must approach the point $z_{dfe}$  when $t\to \I$.


\item[(C2)]  Other conditions are related  to an ``admissible splitting" as a difference of two  parts $\mF,\mV$, called \resp\  ``new infections", and  ``transitions"
\beD A splitting
$$\vi'(t)=\mF(\vi(t),\vec z(t))-\mV(\vi(t),\vec z(t))$$
will be called admissible if
$\mF,\mV$  \saty\ the  following  conditions  \cite{Van08,Shuai}:
\be{cond} \bc \mathcal{F}(\vec 0,\vec z(t))=  \mathcal{V}(\vec 0,\vec z(t))= 0,\\
\mathcal{F}(\vi(t),\vec z(t))\geq 0, \quad \for (\vi(t),\vec z(t)), \\
\mathcal{V}_j(\vi(t),\vec z(t))\le 0, ~~ \text{ when } \vi_j=0,\\
\sum_{j=1}^n \mathcal{V}_j(\vi(t),\vec z(t))\geq 0, \for (\vi(t),\vec z(t)),
\ec\ee
where the subscript $j$ refers the $j$'th component.
\QEDB
\eeD

\beR The splitting of the infectious equations  has its origins in
epidemiology. Mathematically, it is related to the ``splitting of
Metzler matrices"-- see \fe\ \cite{fall2007epidemiological}.
Note however  that the splitting conditions may  be satisfied for several or for no subset of \com s
(see \fe\ the SEIT model, discussed in \cite{Van08}, \cite[Ch 5]{Mart}). Unfortunately, for deterministic epidemic models,  there is no   clear-cut definition of $\mR$ \cite{roberts2007pluses,li2011failure,Thieme}.\fn[3]{A possible explanation is that several stochastic epidemiological  models may correspond in the limit
to the same deterministic  model.}
\QEDB
\eeR

\item[(C3)] We turn now to the last conditions, which concern the linearization of the infectious equations around the DFE. Putting
 $L =A_{i,i}(\vz, \vec z_{dfe})$, and letting $f$ denote the perturbation from the linearization, \wmw:
\beq \la{dynl} \vi'(t) &&=\vi(t) L-f(\vi(t),\vec z(t))=\vi(t)(F-V)-
f(\vi(t),\vec z(t)),\\&&  F:=\left[\frac{\partial \mathcal{F}}{\partial \i}\right]_{x_{dfe}}
 V=\left[\frac{\partial \mathcal{V}}{\partial \i}\right]_{x_{dfe}}, L=F-V. \no \eeq

The ``transmission and transition" linearization matrices  $F,V$
must \saty\ that $F\ge 0$ componentwise and  $V$ is a non-singular M-matrix, which ensures that $V^{-1}\ge 0$.\fn[6]{The  assumption (B) implies that   $L=F-V$ is a ``stability (non-singular) M-matrix", which is necessary for the non-negativity and boundedness of the solutions \cite[Thm. 1-3]{de2019some}.}

 \EEN

Under conditions (C1)-(C3), the \ngm\ method gives an explicit
expression for the basic reproduction number, given by $\mR:=\l_{PF}(FV^{-1})$.

The basic reproduction number is a threshold parameter in the
following sense  \cite[Thm 1]{Van08}:
 \BEN
\item[1(a)] When $\mR<1$, the DFE is locally asymptotically stable;
\item[1(b)] when $\mR>1$, the DFE is unstable;

\item[2~~~~] The DFE is globally asymptotically stable when $\mR \le 1$,
provided the ``perturbation from linearity" $f=i(F-V) - \mF +\mV$ is non-negative \cite[Thm 2]{Van08}.
\fn[5]{Note that the \cite{guo2006global} strong $\mR$ alternative was only established for a general n-stage-progression, which is a particular case of the model we study below, in which  A is an ``Erlang" upper diagonal matrix. It is natural to expect that the result continues to hold for other  non-singular Metzler matrices}
\EEN

In what follows, we will call the alternative 1(a)-1(b)  the
``weak $\mR$ alternative".  In contrast, the result  2 has been called the ``strong $\mR$ alternative" in \cite{guo2006global,Shuai}.

\section{Matrix SIRS epidemics with  demography, loss of immunity, vaccination and
one susceptible class 
 (SIRS epidemics with \PH ``disease time")
\la{s:Feng}}

While the elusive $\mR$  can  be defined as the spectral radius of a certain matrix, provided that the \ngm\ assumptions apply,   it is often possible and  considerably more convenient to employ models where $\mR$ may be
explicitly expressed in terms of the matrices that define
the  model  \cite{ma2006generality,Arino,Feng,Andr}.

The idea behind these models is to further divide  the noninfected \com s   into  S(usceptible)/input classes, defined by producing ``new non-linear infections", and  output R classes (like $D,D_e$ in our example), which are fully determined by the rest, and may therefore be omitted from the dynamics.
 \Fr it is convenient to   restrict to   epidemic models with linear force of  infection, since it is known that non-linear forces of infection may lead to very complex dynamics \cite{liu1986influence,liu1987dynamical,georgescu2007global,tang2008coexistence}, which are not always easy
to interpret epidemiologically.  This is in contrast with the Arino models already studied
  where
  one may typically establish the absence of periodic solutions (closed orbits, homoclinic loops and oriented phase polygons) \cite{Raz,Arinovar}.

 {Below, we study the equilibrium points and the dynamical behavior of such models,  restricting to the case of {\bf one susceptible class}. Our  goal is to understand the effects of
 demography, vaccination and loss of immunity, which are missing in the original paper \cite{Arino}.

\beD  A matrix ``SIRS epidemic" {of type $(1,n,p)$}, with demography parameters $(\L,\mu)$ (scalars), loss of immunity column vector $ \bg_r$ and vaccination row vector  $\vec \g_s$,   is  characterized  by a set of parameters $( A,B, W,  \bn,  {\bn_r})$, where
 $\bn , {\bn_r}$ (in boldface) denote  column  vectors of extra death rates, and $A,B,W$ denote \resp\ matrices of dimensions $n \times n$,  $n \times n$ and $n \times p$, \resp. This model  contains  one  {  susceptible class} $S$, a $p$-dimensional   vector of removed states $\vec R$ (healthy, dead, vaccinated, etc), and a   $n$-dimensional vector of ``disease" states $\vec I $ (which may contain latent/exposed, infective, asymptomatic, etc). {The dynamics are}:

\begin{align}
S'(t)&=\L N  -\fr{S(t)}N  \vec I (t)  \bb -(\vec \g_s+ \mu) S(t) + \vec R(t) \bg_r, \quad\quad \bb= B \bff1,\nonumber\\ 
\vec I'(t)&=  \vec I (t) \pr{\fr{S(t)}N   B + A   -Diag(\bn+ \mu \bff 1)},\la{SYRN}
\\ \vec R'(t)&= \vec  I(t) W + \vg_s S(t) -\vec R(t) (Diag(\bg_r {+\bn_r}+\mu \bff 1)),\nonumber
\\N'(t)&= \S'(t)+ \vec  I'(t) \bff 1 + \vec  R'(t) \bff 1= (\L -\mu) N-
\vec I(t) \bn  {- \vec R(t) \bn_r}, \nonumber
\\D'(t)&=   \mu(  S(t) + \vec I(t) \bff 1 + \vec R(t) \bff 1) ,\nonumber
\\D_e'(t)&=   \vec I(t) \bn  {+ \vec R(t) \bn_r}.\nonumber
\end{align}
In short, we will refer to the above model to as matrix SIR/V+S model.
\QEDB
\eeD

Here, 
\BEN
\im $\vec I(t) \in \mathbb{R}^n$ is a row vector whose components model a set of
disease states (or classes).

\im $\vec R(t) \in \mathbb{R}^p $ is a row vector whose components model a set of
recovered states (or classes), each accounting for individuals who recovered from the
infection.   {In what follows, we will focus on  the case of one recovered class.}

\im $ B $ is a $n \times n$ matrix, where each  {entry} $B_{i,j}$ represents the
force of infection of the  {disease} class $i$ onto class $j$.
 {We will denote by $\bb$ the vector containing the sum of the entries in each row
of $B$, namely, $\bb= B \bff1$.}

\im  $A$ is a $n\times n$ Markovian sub-generator matrix (i.e., a Markovian generator
matrix for which the sum of at least one row is strictly negative),  {where each
off-diagonal entry $A_{i,j}$, $i\neq j$, satisfies $A_{i,j}\geq 0$ and describes the
rate of transition from disease class $i$ to disease class $j$; while each diagonal
entry $A_{i,i}$ satisfies $A_{i,i} \leq 0$ and describes the rate at which individuals
in the disease class $i$ leave towards non-infectious compartments.}
Alternatively, $-A$ is a  non-singular M-matrix \cite{Arino,Riano}.\fn[4]{An M-matrix
is  a real matrix $V$ with  $ v_{ij} \leq 0, \forall i \neq j,$ and having eigenvalues
whose real parts are nonnegative \cite{plemmons1977m}.}

\im $\bn \in \mathbb{R}^n, \bn_r  \in \mathbb{R}^p$  are column vectors describing the death rates in the  {disease and recovered \com s  caused by the epidemic (and possibly vaccinations), \resp.}

\im $\bg_r$ is a vector describing  {the rates at which individuals lose immunity  (i.e. transition from recovered states to the susceptible state).}

\im $\vg_s$ is a vector  {describing the rates at which individuals are vaccinated (immunized).}

\im  $W $ is a ${n \times p}$ matrix  {whose entries model the rates at which
individuals in the disease states transfer to recovered or dead states. In what follows, we assume that the
$n \times (n+p)$ matrix $\T A =\bep A &W \eep $ satisfies $\tilde A \bff 1 = 0$
(namely, the sum of the entries in each row is equal to $0$), which implies mass
conservation.}
\EEN

\beR We have not found any work in the literature on models with linear birth-rate, at this level of generality. As mentioned in the introduction, the immense majority of the literature is dedicated  to
models that may be formally obtained from \eqr{SYRN}  by letting $N=1$ (the idea being that $N$ is approximately constant, either since it is huge, or since it is observed only over a short period of time). We will call this formally obtained model ``classic/pedagogical" where we added the last qualifier to emphasize that it is unrelated to the model we study. This is in contrast with the FA and IA approximations introduced later, which do approximate the scaled version of \eqr{SYRN} introduced below.
\QEDB\eeR

\beR \BEN \im Note that when
$p=1$, then $W=\bff a:=(-A) \bff 1$ is a vector with a \wk\   probabilistic interpretation in the theory of \PH distributions: it is the column vector which completes a matrix with negative row sums to a matrix with  zero row sums.

\im A  particular but revealing case
is that when $p=1$ and matrix $B$ has rank 1, necessarily hence of the form $B=\bb \va$, where $\va\ \in \mathbb{R}^n $ is a {\bf probability row vector} whose components $\a_j$ represent the fractions of susceptibles entering into the  disease compartment $j$, when  infection occurs.  We will call this  SIR-PH, following Riano \cite{Riano}, who emphasized its probabilistic interpretation -- see also \cite{Hurtado}, and see \cite{hyman1999differential} for an early appearance of such models.
\QEDB
\EEN
\eeR

\beR  Note in the general model the factorization of the equation for the diseased \com s, which turns out to be an essential
feature of the model.

 Note also that our model includes  important epidemiological  parameters,  {such
as $\bn$ (describing the death rate when individuals are in the infectious
compartments) and $\bg_r$ (describing the rates at which recovered individuals lose
immunity), which are often omitted in simpler models. In what follows, our purpose  is
precisely to study the emerging behaviors due to the presence of these parameters with
respect to the simple Arino model.}
\QEDB
\eeR

It is  convenient to reformulate \eqr{SYRN} in terms of the fractions
normalized by the total population
\begin{align}
\label{eq:fractions}
\s&=\fr{S}N, &
\vi&=\fr{1}N \vec I,&
\vr&=\fr{1}N \vec R N= \s + \vi . \bff 1 + \vr.\bff 1 .
\end{align}

The reader may check that  the following equations hold for  the scaled variables:
\be{SYRsc}
\bc
 \s'(t)= \L -\pr{ \L + \g_s } \s(t) + \vr(t) \bg_r- \s(t)\; \vi(t) \pr{ \bb -  \bn}   { + \s(t) \vr(t) \bn_r}\\
  \vi '(t)=    \vi (t) \Big(\s(t)  \; B +\pr{\vi(t) \bn+  {\vr(t) \bn_r}}I_n+ A- {Diag\pp{\bn+\L}}
 \Big)  \\
\vr'(t) =   \s(t) \vec \g_s +  \vi (t) W  - \vr(t) \Big( {Diag \pp{\bg_r+\bn_r +\L} -(\vi(t) \bn+ \vr(t) \bn_r) }I_p\Big)\\
\s(t) + \vi(t) \bff 1 + \vr(t) \bff 1=1
\ec,
\ee
where we let $\g_s:=\vec \g_s \bff 1$.
Moreover, by letting $\n:=\s+\vi \bff 1+\vr \bff 1,$ we have
$$\n'(t)=(\L -\vi(t) \bn)(1-\n(t))=0.$$
Hence, the above equation guarantees that  if
$\s(t_0)+\i(t_0)+\r(t_0)=1$ for some $t_0 \in \realnneg$, then
$\s(t)+\i(t)+\r(t)=1$ for all $t \geq t_0$.
Accordingly, in what follows we will always assume that $n(t_0) =1$,
which guarantees that $\n(t)=1, \for t$.

The following definition puts in a common framework the dynamics for the scaled
process and two interesting approximations.

\beD \la{d:fisg} Let $\Phi_s, \Phi_i, \Phi_r \in \{0,1\}$ and let
\begin{align}
\label{SYRsc-def}
 \s'(t) &= \L -\pr{ \L + \g_s } \s(t) + \vr(t) \bg_r
- \s(t) \vi(t) \bb
+ \Phi_s \s(t) ( \vi(t) \bn    +  \vr(t) \bn_r),
\nonumber\\
\vi '(t) &=  \vi (t) \Big(\s(t)  \; B
+ A
- {Diag\pp{\bn+\L \bff 1}}
\Big)
+ \Phi_i \vi(t) \left(\vi(t) \bn+  {\vr(t) \bn_r} \right),
\nonumber\\
\vr'(t) &=   \s(t) \vec \g_s +  \vi (t) W  - \vr(t)  {Diag \pp{\bg_r+\bn_r +\L \bff 1}  }
+ \Phi_r \vr(t) (\vi(t) \bn+ \vr(t) \bn_r),
\nonumber\\
\s(t) &+ \vi(t) \bff 1 + \vr(t) \bff 1=1.
\end{align}

\BEN
\im The model \eqref{SYRsc-def} with $\Phi_s=\Phi_i=\Phi_r=1$  will
be  called scaled model (SM).

\im  The model \eqref{SYRsc-def} with $\Phi_s=\Phi_i=\Phi_r=0$  will
be called first approximation (FA).

\im  The model \eqref{SYRsc-def} with $\Phi_s=\Phi_r=1$ and $\Phi_i=0$
will  be called  intermediate approximation (IA).
\QEDB
\EEN
\eeD

Fig.~\ref{f:PIP} compares the qualitative behavior and equilibrium points of the $(\s,\i)$ coordinates of the three variants
of a SIR-type example of \eqr{SYRN}  model (discussed  in detail in \cite{AABGH}).

\begin{figure}[H]
\centering
\includegraphics[width=.8\columnwidth]{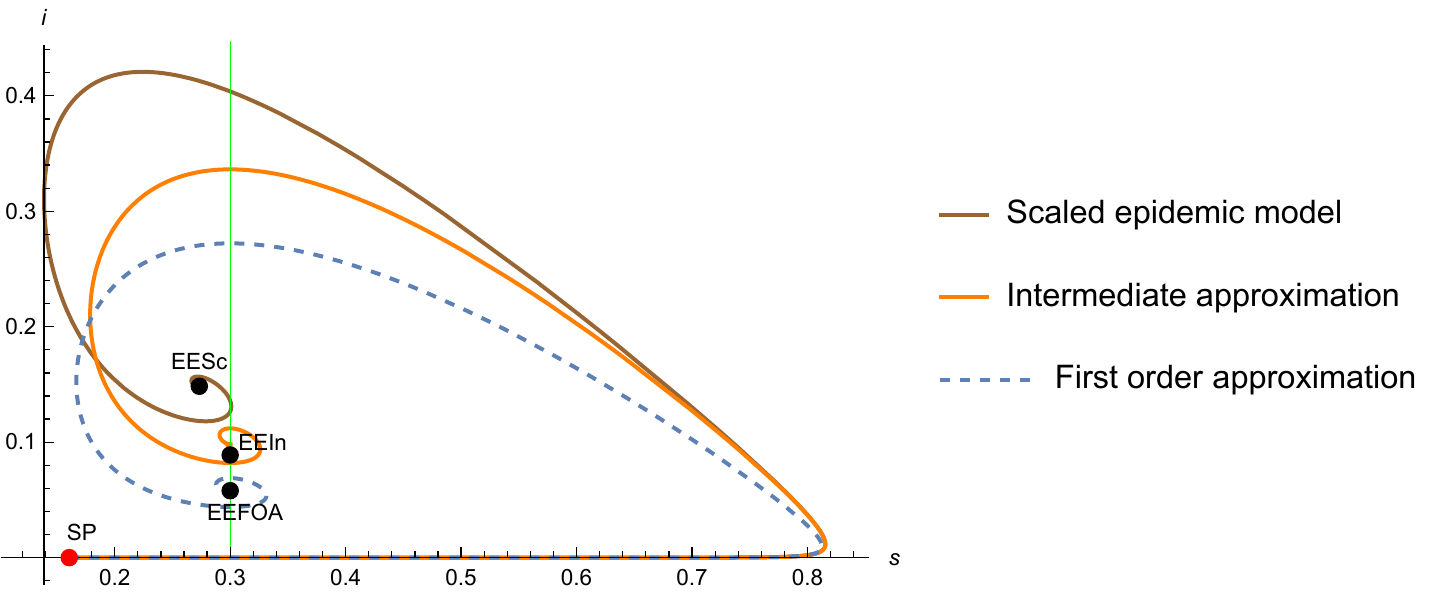}
\caption{Parametric $(\s,\i)$ plots of the scaled epidemic and   its FA and intermediate approximations for a SIR-type model with one infectious class, starting from a starting point SP with $i_0=10^{-6}$, with
 $\mR=3.21, \b=5,\; \g  = 1/2,\; \L=\mu=1/10, \; \g_r=1/6, \; \g_s=.01, \;\nu_i=.9, \; \nu_r=0$. $EESc, EEIn,EEF0A$ are the stable endemic points of the scaled model, intermediate model, and the FA model, respectively.
 The green vertical line denotes the immunity threshold $1/\nR= s_{EEF0A}= s_{EEIn}.$
 Note that    the epidemic will  spend at first a long time (since births and deaths have  slow rates as compared to the disease) in  the vicinity of the manifold $\vi(t) =0$, where the three processes are indistinguishable, before turning towards the  endemic  equilibrium point(s).}
 \label{f:PIP}
\end{figure}

\begin{example}  The classic SEIRS  model
 is a particular case of SIR-PH-FA model (see also next remark)
 obtained when $$\va=(1,0),A=\bep -\g_e&\g_e\\0&-\g \eep, W=\bff a
 =\bep 0\\\g \eep, \bb=\bep 0\\\b\eep, \; \mbox{so}\; B=\bep 0 & 0\\  \b&0 \eep , \bn=\bep 0\\\nu\eep;$$ this yields: \be{seirs} \bc
\s'(t)= \L   -\s(t) \pr{\b \i(t)+\g_s+{\L}} + \g_r \r(t)\\
\bep \e'(t)\\\i'(t)\eep =
\bep  {-(\g_e +\L)} &  {\b \s(t)} \\
 \g_e& -\pr{\g   +\L +\nu}\eep
  \bep \e (t)\\\i(t)\eep
\\
\r'(t)=  \g_s \s(t)+ \g   \i(t)- \r(t)(\g_r +\L +\nu_r ) \ec .
\ee

\begin{figure}[H]
\centering
\includegraphics[scale=0.7]{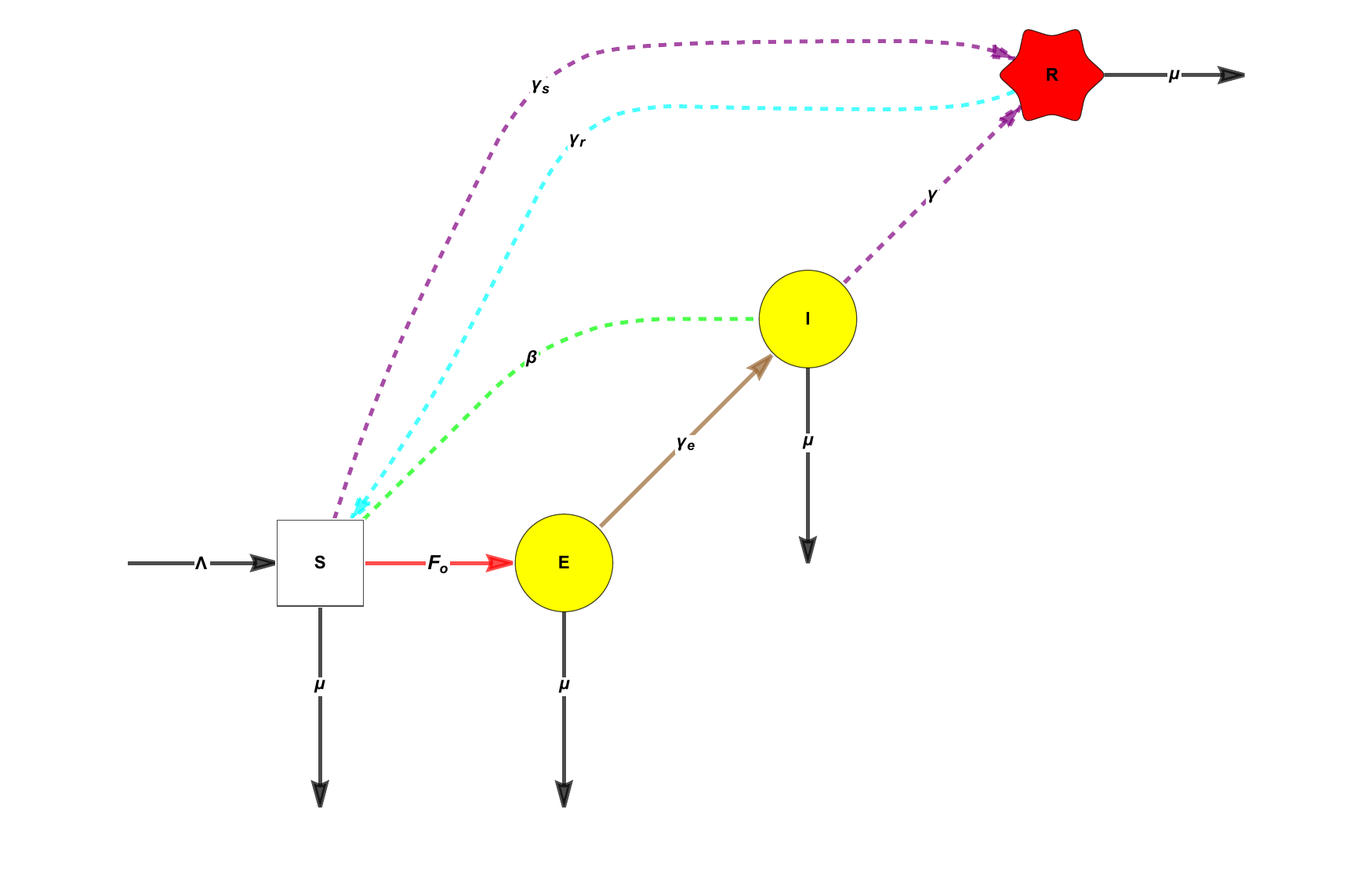}
\caption{Chart flow of the SEIRS model. The forces of infection  is $ F_o=\b \i.$ The red edge corresponds to the entrance of susceptibles into the disease classes, the dashed green edges correspond to contacts between  diseased to susceptibles, the brown edge is the rate of the transition matrix $V$, and the cyan dashed line corresponds to the rate of loss of immunity. The remaining black lines correspond to the inputs and outputs of the birth and natural death rates, respectively. \label{f:SEIRS}}
\end{figure}

\end{example}

\beR From now on we will write the dynamical systems with the vector state variables pre-multiplying the model matrices.
The reason is that we want   $\va$  in the SIR-PH fundamental case to be a row vector, as in the theory of \PH\ distributions, which turns out to be convenient in applications.
\QEDB
\eeR

 \beR  We point out now that the {scaled matrix SIR/V+S model} belongs to a class of  models introduced, for different motivations, in \cite{fall2007epidemiological}. Indeed, after dropping  the first equation $$\s'(t)=\L-\s(t) \vi(t) \bb -(\L+\g_s )\s(t)+\vr(t) \bg_r +\pr{\vi(t) \bn   + \vr(t) \bn_r} \s(t)$$ and rewriting the rest as:
\bea &&\bep \vi'(t)&\vr'(t)\eep = \bep 0 & \s(t) \bg_r    \eep+ \bep\vi(t)&\vr(t)\eep
\pp{\bep
\s(t)  B+A -  (\L +\nu) I_n&W  \\
  \bff 0 & {-Diag(\bg_r)+\L I_p}\eep -Diag(\bff \nu_x)} \\&&+\vec x Diag(\bff \nu_x) Diag( \vec x), \eea
  where $\vec x=(\vi, \vr), \bff \nu_x=(\bff \nu, \bff \nu_r)$, we recognize a particular case of \cite[(2.1)]{fall2007epidemiological}.
Finally, we note that the {\bf FA} is obtained by dropping the row of the last matrices.
\QEDB
\eeR


Using now $\nabla \vi(t) \bn \vi(t)= \bn \vi(t)+\vi(t) \bn I_n$, and putting $x=\vi \bn+ \vr \bn_r$, we find that the transpose of the Jacobian matrix of the scaled model \eqr{SYRsc} is given by
{\small
\bea J=
\bep  x-\vi \bb-(\L+\g_s) &\vi B   &  \vec \g_s\\
\s\pr{\bn- \bb} & \s B  + A-Diag(\bn+ \L \bff  1)+\bn \vi +(\vi \bn +\vr \bn_r)I_n &W + \bn \vr\\
\g_r  {+ s \bn_r} & \bn_r \vi & { \bn_r \vr-Diag[\bg_r +\bn_r+\L]+(\vi \bn+\vr \bn_r )I_p}\eep.
\eea
}
 {At the disease-free equilibrium $(s_{dfe}, \vec 0, \vr_{dfe})$, we have}

\bea
J_{dfe}=
 \bep  -(\L+\g_s)  {+\vr_{dfe} \bn_r}& 0 & \vec \g_s \\ s_{dfe} \pr{\bn- \bb} & s_{dfe} B  + A-Diag(\bn+ \L \bff  1) +\vr_{dfe} \bn_rI_n& W+ \bn \vr_{dfe}\\   {\g_r + s_{dfe} \bn_r} &0 & { \bn_r \vr_{dfe}-Diag[\bg_r +\bn_r+\L]+(\vr_{dfe} \bn_r )I_p} \eep.
\eea

\beR a) This is (up to a transpose) a generalization of the Jacobian of the scalar SIR/V+S model, with the middle element replaced by a $n\times n$ matrix.

b) Note the structure of the second column, which is equivalent to the factorization property, and implies
Proposition \ref{p:brn}, below. \QEDB\eeR


\section{Stability results for the SIR-PH model \la{s:S}}

We consider here  the particular case  with a single
recovered class, called  SIR-PH model, where $W=\bff a =(-A) \bff 1$.
Notice that, in this case, we have
a reduced set of parameters $(\L=\mu,  \va, A, \bb, \bn, {\nu_r}, \g_s, \g_r)$.

\ssec{The  \brn\ for  SIR-PH, via the \ngm\ method \cite{Van,Diek}
\la{s:R0}}

We follow up here on a remark preceding \cite[Thm 2.1]{Arino}, and show in the following proposition that their simplified formula for the \brn\ still  holds when loss of immunity and vaccination are allowed, provided that $p=1$ and $B=\bb \va$ has rank one.

\beP \la{p:brn}  Consider a SIR-PH model (i.e. a single recovered class),
 with parameters
$(\L=\mu,  \va, A, \bb, \bn,  {\nu_r}, \g_s, \g_r)$, and matrix of recovery rates $W=\bff a =(-A) \bff 1$.

\BEN \im
 \BEN \im When $\nu_r=0$,  the unique  DFE is $(\fr{\L +\g_r}{\L+\g_r+\g_s},0,\fr{\g_s}{\L+\g_r+\g_s})$.


   \im When $\nu_r>0$, exclude the case $\L=\g_r=0$. Then, there exists a unique  DFE $(\sd,0,1-\sd)\in \mD$, where $\sd$ \sats\
 the second order equation
$\s \pp{ \nu_r  \s  +\L+\g_r+\g_s-\nu_r}-(\L +\g_r)=0 $ and is given by
\be{defsc}
 s_{dfe}=  \frac
{\sqrt{4 \nu _r \left(\Lambda +\gamma _r\right)+
\left(\Lambda +\gamma _r+\gamma
   _s-\nu _r\right)^2} -
   \pr{\Lambda +\gamma _r+\gamma _s-\nu _r}}
   {2 \nu _r}.
\ee

\EEN
\im   The weak $\mR$ alternative holds for
the threshold parameter
\be{stab} \mR= \l_{PF}(F V^{-1}),\ee
 where $F,V$ are defined in \eqr{FV},  and $\l_{PF}$ denotes the (dominant) \PF\ eigenvalue.

\im  For rank one, $B:=\bff b \va$  and $\nu_r=0$, we further have
\BEN \im
\be{nRV}\mR= s_{dfe}\; \nR, \text{ where } \nR=   \va\ V^{-1} \; \bb.   \ee

\im  If $\mR \leq 1$, and if the perturbation from linearity
 defined in \eqr{dynl} is \nne, then the scalar combination  \be{Y}Y=   \vi \;  V^{-1} \; \bb   \ee
 is a Lyapunov function for the DFE.
 \EEN

 \EEN
\eeP

\Prf 1. The disease free system ( {with $\i=0,\r=1-\s$}) reduces to
\be{SIRDF}
\s'(t)= \L   -(\g_s+\L) \s(t) +( \g_r+\nu_r \s(t)) (1-\s(t)).
\ee
For the fixed points we must, depending on $\nu_r$,  solve either a quadratic, or a linear  equation
$$\bc \L +\g_r- \s \pp{ \nu_r  \s  +\L+\g_r+\g_s-\nu_r}=0& \nu_r >0\\
\s \pp{\L+\g_r+\g_s}-(\L +\g_r)=0& \nu_r =0\ec.$$

 One root
\be{sdfe}
 s_{dfe}= \bc \fr{\L +\g_r}{\L+\g_r+\g_s}& \nu_r =0\\ \frac
{\sqrt{\D_{dfe}} -
   \pr{\Lambda +\g_r+\g_s-\nu _r}}
   {2 \nu _r}, \qu \D_{dfe}=4 \nu _r \left(\Lambda +\g_r\right)+
\left(\Lambda +\g_r+\gamma
   _s-\nu _r\right)^2&\nu_r > 0 \ec
\ee
is always  in $[0,1]$ and will be denoted  by $s_{dfe}$.  

\beR $s_{dfe}$ is continuous in $\nu_r$, since for $\nu_r$ small,
$s_{dfe} \approx \frac
{\Lambda +\g_r+\gamma
   _s-\nu _r+\fr{2 \nu _r \left(\Lambda +\g_r\right)}{\Lambda +\g_r+\gamma
   _s-\nu _r} -
   \pr{\Lambda +\g_r+\g_s-\nu _r}}
   {2 \nu _r} \to \fr{\Lambda +\g_r}{\Lambda +\g_r+\gamma
   _s}$
   (this approximation may be made rigorous by applying the rule of l'Hospital).
\QEDB
\eeR

 The other root in the  quadratic case $\nu_r>0$
\be{dfe2}
  \frac
{\nu _r- \pr{\Lambda +\g_r+\g_s}-\sqrt{4 \nu _r \left(\Lambda +\g_r\right)+
\left(\Lambda +\g_r+\gamma
   _s-\nu _r\right)^2}
  }
   {2 \nu _r}
\ee
is strictly negative, unless
\be{bop}\bc \L+\g_r=0\\ \nu_r \geq \g_s + \L+\g_r\ec \Eq \bc \L=\g_r=0\\  \nu_r \geq \g_s\ec,\ee in which case it yields  a second DFE point with $\s=0=\i$, which we exclude (in order to be able to  apply the \ngm\ method).

2.  It is enough to show here that the conditions of \cite[Thm 2]{Van08} hold, with respect to the infectious set $\vi$, and a certain splitting.

The DFE and its local stability for the disease-free system have already been checked in the SIR/V+S example.


We provide now a splitting for the infectious equations:
\bea \vi'(t) &=& \vi(t)\pp{\s(t)  B + \vi(t) \bn I_n + \nu_r \r(t) I_n}
-\vi(t)\pp{Diag(\bn+ \L \bff  1)-A }:=\mathcal{F}(s,\vi)-
\mathcal{V}(\vi) \eea
 (where $\r=1-\s- \vi \bff 1$).
 The  \corr\ gradients at the DFE $\vi=0$ are
\be{FV} \bc F= \left[\frac{\partial \mathcal{F}(X^{(DFE)})}{\partial \vi}\right]=\s B +\nu_r\r \; I_n,\\
 V=  \left[\frac{\partial \mathcal{V}(X^{(DFE)})}{\partial \vi}\right]= Diag(\bn+\L \bff 1)-A.
 \ec \ee
We note that $F$  has  non-negative \elt s, and that $V$ is a  M-matrix, and therefore $V^{-1}$  exists and has  non-negative \elt s, $\for \L, \bn$.
We may check that the conditions  \eqr{cond} are \satd.

For example,  the  last non-negativity condition in \eqr{cond}
 \be{par} \vi(t)\pp{Diag(\bn+ \L \bff  1)-A  } \bff 1 \ge 0, \for \vi \in \mD, \ee
is  a consequence of $-A$ being a  M-matrix, which implies $-A \bff 1 \ge 0,$ componentwise.

3.a)   Now  if $n=1$,  or if $B= \bb \va$ has rank 1, and $\nu_r=0$, the matrix $F$ in \eqr{FV} is the product of a column vector and a row vector, the dimension of  its image is one, and the same holds for $B V^{-1}$. Equivalently, $rank(B V^{-1})=1$.
The \textit{"rank-nullity theorem"} $rank(B V^{-1})+ nullity(B V^{-1})= n$
\cite{horn2012matrix} implies that $(n-1)$ of the eigenvalues of $B V^{-1}$ are zero, and the Perron-Frobenius eigenvalue is the remaining one.
This latter one must be equal  to the trace of $ B V^{-1}= \bb \; \va\ V^{-1}$ which may be checked to equal   $ \va\  V^{-1}\; \bb$.  Finally, the linearity $\mR= s_{dfe}\; \nR $  is obvious.

3.b) is a particular case of \cite{Shuai}, since the \PF\ eigenvector  in our rank-one case may be taken as $\bff b$.
\QED

\ssec{The classic/pedagogical  SIRS-FA model \label{s:EqP}}

In this section, we give more explicit results for the \DFE\ and the endemic
equilibrium of th following model, referred to as SIRS-FA,
\be{SYRscF}
\bc
 \s'(t)= \L -\pr{ \mu + \g_s } \s(t) + \vr(t) \bg_r- \s(t)\; \vi(t)  \bb   \\
  \vi '(t)=    \vi (t) \Big(\s(t)  \; \bb \va - V\Big)  \\
\vr'(t) =   \s(t) \vec \g_s +  \vi (t) W  -
\vr(t) \Big( Diag \pp{\bg_r+\bn_r +\mu \bff 1} \Big)\\
\s(t) + \vi(t) \bff 1 + \vr(t) \bff 1=1
\ec,
\ee

\beP \la{p:ee}  \BEN \im The { pedagogical $(\L, \mu, A, \vg_s,\bg_r)$} SIRS  system \eqr{SYRscF} has a unique disease-free equilibrium (DFE) fixed point
 \bea(s_{dfe},\vz, \vr_{dfe}), \quad  s_{dfe}=\frac{\L}{\mu +\g_s-\vg_s \pr{Diag({\bg_r+\bn_r + \mu \bff 1})}^{-1}  \bg_r}, \quad \vr_{dfe}=s_{dfe} \vg_s Diag({\bg_r+\bn_r +\mu \bff 1})^{-1}.\eea

 In the SIR-PH case, the DFE simplifies to\fn[4]{This formula has appeared already in many particular cases --see \fe\ \cite[(19-20)]{Sen}}: \be{sr} (s_{dfe},\vz, r_{dfe})=\left(\frac{\L(\mu+\g_r+\nu_r)}{\mu \g_r+(\mu+\nu_r)(\mu+\g_s)}, \vz, s_{dfe} \g_s \rt).
 \ee

 \im  If $\mR >1$, then the { pedagogical   system} \eqr{SYRscF} has a unique second fixed point within its forward-invariant set. 
This endemic fixed point  is  such that
$1/s_{ee}$ is an eigenvalue of the matrix $B V^{-1}$.

 In the SIR-PH case  it must \saty\
 $$ 1/s_{ee}={\nR}=\va\ V^{-1}\mathbf{b}.$$

The disease components $ \vi_{ee}$ \saty:
\be{ino}  \vi_{ee} (\fr 1{\nR}  \; B  + A- Diag(\bn +\mu \bff  1)):= \vi_{ee} M=0   \ee
 ( $\vi_{ee}$ is a Perron-Frobenius eigenvector of the  matrix $M$ related to the \ngm).

\im  The normalization of $i_{ee}$ is given by \eqr{ieen} below.
When $ \mu=\L $, 
this becomes:
\be{ino}   \vi \pr{\bb- \nR \bff a \frac{ \g_r}{\g_r+\nu_r+\L}}= \L \nR -\L +\g_s \pr{\frac{\g_r}{\L+ \g_r+\nu_r} -1}.
 \ee

\im The disease-free equilibrium is locally asymptotically stable if $\mR< 1$ and is unstable if $\mR > 1$.

\im When $\nu_r=0$, the critical vaccination   defined by solving $\mathcal{R}_0=1$ with respect to $\g_s$ is given by
\be{crv}  \g_s^*:= (\L+\g_r)\pr{\va V^{-1} \bff b -1}= (\L+\g_r)\pr{\nR -1}. \ee

\EEN
\eeP

We show in Fig.~\ref{f:spFOA} a stream plot of the SIRS-FA model that illustrates the above properties.

\begin{figure}[H]
\centering
\includegraphics[scale=0.9]{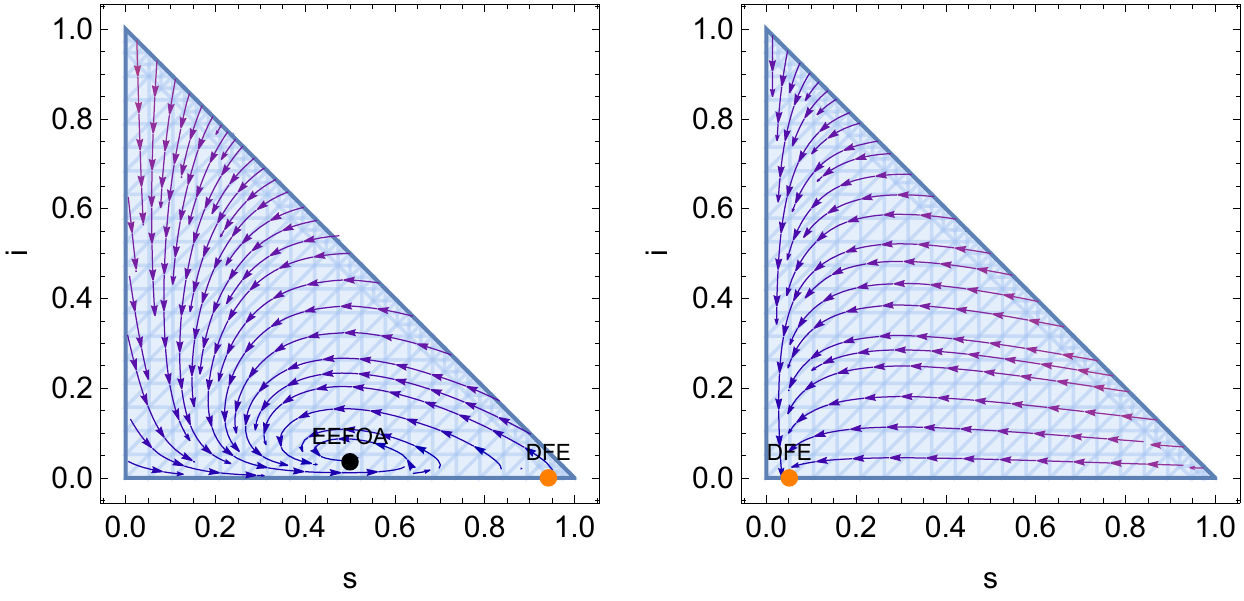}
\caption{Stream plots of $(s,i)$ for an example of  SIRS-FA model with one infectious class \cite{AABGH},    when $\g_s \in \left\{1/100,  3\right\}$ is smaller and bigger, \resp, than the critical vaccination $\g_s^*=0.239087$ defined in \eqr{crv}.}
\label{f:spFOA}
\end{figure}

\Prf The proof starts by examining the two cases  which arise from factoring the disease equations. \Mp  we will search separately in the {\bf disease free set} $\mI:=\{\vi=\vz\}$ and in its complement $\mI^c$. Then:
\BEN \im  either $\vi\in \mI$ and  solving
\bea
\bc
 \L = \pr{ \mu + \g_s } \s - \vr \bg_r  \\
\vz =    \s \vec \g_s - \vr Diag(\bg_r+\bn_r + \mu \bff  1)
\ec=\bep \s &\vr \eep \bep { \mu + \g_s } & \vec \g_s  \\
-  \bg_r&  - Diag(\bg_r+ \bn_r + \mu \bff  1)\eep
\eea
for $\s, \vr$ yields the unique DFE  (it may be shown by induction that the determinant is negative). Or,

\im  the determinant of the resulting homogeneous linear system  for $\vi \neq \vz$ must be $0$, which implies that $\s=s_{ee}$ \sats\ \be{lP} det\pp{s_{ee}  \; B  + A -Diag(\bn +\mu \bff  1)}=0. \ee

Note now that $V=Diag(\bn +\mu \bff  1)-A$ is an invertible matrix.
Using   $det (U U') = det( U)  det( U')$, \eqr{lP} \mbw \bea  det \left[\pr{s_{ee} B - V} V^{-1}\right] = 0 \Eq det \left[s_{ee}B V^{-1} - \mathbf{I}\right] = 0\eea

 Dividing then  by $s_{ee}$ yields the characteristic polynomial of a matrix
\be{ls}
det \left[B  V^{-1} -\frac{1}{s_{ee}} \mathbf{I}\right] = 0.
\ee

In the SIR-PH case, noting that the rank one matrix $B  V^{-1}$ has  $(n-1)$  eigenvalues equal to zero, we conclude that
 the inverse of the susceptible fraction  $1/s_{ee}$ of an endemic state must equal the \PF\ eigenvalue $\nR$.
Note that $s_{ee} <1$   follows from our assumption on $\nR$. The other coordinates are determined starting with $\vi$, which must be proportional to a Perron-Frobenius nonnegative eigenvector.

 \beR {Equivalently, \itc\ \eqr{ls} may be reformulated as saying that $s_{ee}$ must \saty
$$\l_{PF}(s_{ee}  \; B  \bep A -Diag(\bn +\mu \bff  1))\eep^{-1} =1 \Eq \l_{PF}(s_{ee}  \; B  + A -Diag(\bn +\mu \bff  1))=0,$$ where $\l_{PF}$  denotes the Perron-Frobenius eigenvalue. We note that the matrices in the last formulation intervene also in the \ngm\ approach. }
\QEDB
\eeR

 \im   Recall the system
\bea
&&\bc
 0= \L -\pr{ \mu + \g_s } \s + \vr(t) \bg_r- \s\; y  \\
  0=    \vi  (\s  \; B  + A- {Diag(\bn+ \mu \bff  1)}) \\
0 =  \vi  W +  \s \vec \g_s - \vr (Diag(\bg_r+\nu_r+ \mu \bff  1))
\ec. \eea

Since  $s_{ee}=\fr 1{\nR}$, and $\vi_{ee}$ is known up to the proportionality  constant $y=\vi \bb$,  it only remains to solve the last equation in \eqr{ieen} below:
\be{ieen}\bc \vr =  \pp{\vi  W +  \s \vec \g_s } Diag({\bg_r+\nu_r+\mu \bff 1})^{-1}\\
 \vi \pr{\s  \bb-W Diag({\bg_r+\nu_r+\mu \bff 1})^{-1} \bg_r}  = \L -\pr{ \mu + \g_s  } \s +   \s \vec \g_s   Diag({\bg_r+\nu_r+\mu \bff 1})^{-1} \bg_r
 \ec
\ee

This equation may be solved numerically.  When
$p=1 \Lra W:=\bff a=(-A) \bff 1,$  the last formula yields
\be{iees}
\vi \pr{\bb- \nR \bff a \frac{ \g_r}{\g_r+\nu_r+\mu}}= \L \nR -\mu +\g_s \pr{\frac{\g_r}{\g_r+\nu_r+\mu} -1}
\ee

In the particular case $ \mu=\L $, this yields
\eqr{ino}.

\im  This  follows from \cite[Theorem 2.1]{Shuai} (it is a consequence of the fact   that  a linear  function proportional to the  associated Perron eigenvector is a Lyapunov function when $ \mR <1$).


, with\im The result is immediate by solving $\mR=1$ \wrt\ $\g_s$, where $\mR$ is defined in \eqr{nRV}.
\QED
\EEN

\subsection{A glimpse of the intermediate approximation model for matrix SIRS, with $\L=\mu$}

The intermediate approximation associated to \eqr{SYRsc} is
\be{SYRIA}
\bc
 \s'(t)= \L -\pr{ \L + \g_s } \s(t) + \vr(t) \bg_r- \s(t)\; \vi(t) \pr{ \bb -  \bn}   + \s(t) \vr(t) \bn_r\\
  \vi '(t)=    \vi (t) \Big(\s(t)  \; B + A- {Diag\pp{\bn}-\L I_n}
 \Big)  \\
\vr'(t) =   \s(t) \vec \g_s +  \vi (t) W  - \vr(t) \Big(Diag \pp{\bg_r+\bn_r+(\L -(\vi(t) \bn+ \vr(t) \bn_r))\bff 1}\Big)
\ec;
\ee
when $i=0$, we have
\be{siaDFE} s_{dfe}= \frac{\L+\vr_{dfe} \bg_r}{\L+ \g_s -\vr_{dfe} \bn_r} ,\ee
and from the last equation in \eqr{SYRIA}, $\vr_{dfe}$ satisfies the following {third order} equation

\bea
 \pr{\L+ \g_s -\vr_{dfe} \bn_r}\vr_{dfe} \Big(Diag \pp{\bg_r+\bn_r+(\L- \vr_{dfe} \bn_r) \bff 1} \Big)- \vr_{dfe} \bg_r \vec \g_s=\L \vec \g_s
\eea

\beP \la{p:IA} Assume $\bn_r=0$. Then:
 a) The DFE points of the scaled, the intermediate approximation, and the FA are equal,
 with
 $$\vr_{dfe}=\L \vec \g_s \Big((\L+ \g_s)  Diag \pp{\bg_r+\L \bff 1} -  \bg_r \vec \g_s\Big)^{-1}.$$

 b) In the SIR-PH case with $\vec r$ scalar, they reduce all to $(\fr{\L + \g_r}{\L+\g_r+\g_s},0,\fr{ \g_s}{\L+\g_r+\g_s})$ \eqr{sdfe}.

 c) The endemic point is unique. It \sats\ $\s_{ee}=\fr 1 \nR$, $\vi$ is an \eigv\ of the matrix $\fr 1 \nR\; B + A- Diag\pp{\bn}-\L I_n$ for the \eig\  $0$, and
 $$\vr_{ee}=(\s_{ee} \vec \g_s + \vi_{ee} W)\Big( Diag \pp{\bg_r+(\L - \vi_{ee} \bn)\bff 1} \Big)^{-1}$$
 and it \sats\ the normalization
 \be{noI}\vi(t) \pr{ \bb -  \bn} - \nR  \vi W \Big( Diag \pp{\bg_r+(\L - \vi \bn)\bff 1} \Big)^{-1} \bg_r =\L \pr{ \nR-1} - \g_s  + \vg_s \Big( Diag \pp{\bg_r+(\L - \vi \bn)\bff 1} \Big)^{-1} \bg_r  \ee
\eeP
\Prf a) The equations determining the three DFE's coincide.

b) When $\vec r$ is scalar, we find $$\r_{dfe}=\L  \g_s \fr 1{(\L+ \g_s)   \pp{\g_r+\L } -  \g_r  \g_s}= \fr {\g_s }{\L+ \g_s+\g_r }.$$

{c) We have
\bea
\bc
 0= \L -\pr{ \L + \g_s } s_{ee} + \vr_{ee} \bg_r- s_{ee}\; \vi_{ee} \pr{ \bb -  \bn},\\
 \vr_{ee}=(\s_{ee} \vec \g_s + \vi_{ee} W)\Big( Diag \pp{\bg_r+(\L - \vi_{ee} \bn)\bff 1} \Big)^{-1},\\
 \ec
\eea
by susbtitution we get
\bea
\L(\nR-1)-\g_s +\vec \g_s \Big(Diag \pp{\bg_r+(\L - \vi_{ee} \bn)\bff 1} \Big)^{-1}\bg_r=\vi_{ee}(\bb-\bn)-\nR \vi_{ee}W \Big( Diag\pp{\bg_r+(\L - \vi_{ee} \bn)\bff 1} \Big)^{-1}\bg_r.
\eea }
\QED
\beR
{When $p=1$ 
and  $\vi \bn=0$, \eqr{noI} reduces to \eqr{ino} when $\nu_r=0$.}
\eeR

\section{Conclusions and further work}
We have   provided here a few general results for Arino models with varying population, and one susceptible class only. The following directions seem worthy of further work.
\BEN \im The case of two or more susceptible classes.
\im The determination of the largest domain of attraction of the DFE, through which some linear Lyapunov function decreases, might also be approachable via linear programming.
\EEN


\section*{Acknowledgments}

We thank N. Bacaer for providing the references
\cite{bacaer2021mathematiques,lotka1939analyse}, and the referees for their
reviews and suggestions.
%
%
%

\section*{Declarations}

\noindent
\textbf{Competing interests.} The authors declare no affiliation or involvement in any organization or entity with a financial or non-financial interest in the subject matter or materials discussed in this manuscript.

\noindent
\textbf{Funding.}
This work was supported in part by the AB Nexus seed grant from the University of Colorado.





\bibliography{Pare38}


\begin{thebibliography}{67}
\ifx \bisbn   \undefined \def \bisbn  #1{ISBN #1}\fi
\ifx \binits  \undefined \def \binits#1{#1}\fi
\ifx \bauthor  \undefined \def \bauthor#1{#1}\fi
\ifx \batitle  \undefined \def \batitle#1{#1}\fi
\ifx \bjtitle  \undefined \def \bjtitle#1{#1}\fi
\ifx \bvolume  \undefined \def \bvolume#1{\textbf{#1}}\fi
\ifx \byear  \undefined \def \byear#1{#1}\fi
\ifx \bissue  \undefined \def \bissue#1{#1}\fi
\ifx \bfpage  \undefined \def \bfpage#1{#1}\fi
\ifx \blpage  \undefined \def \blpage #1{#1}\fi
\ifx \burl  \undefined \def \burl#1{\textsf{#1}}\fi
\ifx \doiurl  \undefined \def \doiurl#1{\url{https://doi.org/#1}}\fi
\ifx \betal  \undefined \def \betal{\textit{et al.}}\fi
\ifx \binstitute  \undefined \def \binstitute#1{#1}\fi
\ifx \binstitutionaled  \undefined \def \binstitutionaled#1{#1}\fi
\ifx \bctitle  \undefined \def \bctitle#1{#1}\fi
\ifx \beditor  \undefined \def \beditor#1{#1}\fi
\ifx \bpublisher  \undefined \def \bpublisher#1{#1}\fi
\ifx \bbtitle  \undefined \def \bbtitle#1{#1}\fi
\ifx \bedition  \undefined \def \bedition#1{#1}\fi
\ifx \bseriesno  \undefined \def \bseriesno#1{#1}\fi
\ifx \blocation  \undefined \def \blocation#1{#1}\fi
\ifx \bsertitle  \undefined \def \bsertitle#1{#1}\fi
\ifx \bsnm \undefined \def \bsnm#1{#1}\fi
\ifx \bsuffix \undefined \def \bsuffix#1{#1}\fi
\ifx \bparticle \undefined \def \bparticle#1{#1}\fi
\ifx \barticle \undefined \def \barticle#1{#1}\fi
\bibcommenthead
\ifx \bconfdate \undefined \def \bconfdate #1{#1}\fi
\ifx \botherref \undefined \def \botherref #1{#1}\fi
\ifx \url \undefined \def \url#1{\textsf{#1}}\fi
\ifx \bchapter \undefined \def \bchapter#1{#1}\fi
\ifx \bbook \undefined \def \bbook#1{#1}\fi
\ifx \bcomment \undefined \def \bcomment#1{#1}\fi
\ifx \oauthor \undefined \def \oauthor#1{#1}\fi
\ifx \citeauthoryear \undefined \def \citeauthoryear#1{#1}\fi
\ifx \endbibitem  \undefined \def \endbibitem {}\fi
\ifx \bconflocation  \undefined \def \bconflocation#1{#1}\fi
\ifx \arxivurl  \undefined \def \arxivurl#1{\textsf{#1}}\fi
\csname PreBibitemsHook\endcsname

\bibitem{Arino}
\begin{barticle}
\bauthor{\bsnm{Arino}, \binits{J.}},
\bauthor{\bsnm{Brauer}, \binits{F.}},
\bauthor{\bparticle{van~den} \bsnm{Driessche}, \binits{P.}},
\bauthor{\bsnm{Watmough}, \binits{J.}},
\bauthor{\bsnm{Wu}, \binits{J.}}:
\batitle{A final size relation for epidemic models}.
\bjtitle{Mathematical Biosciences \& Engineering}
\bvolume{4}(\bissue{2}),
\bfpage{159}
(\byear{2007})
\end{barticle}
\endbibitem

\bibitem{Van08}
\begin{bchapter}
\bauthor{\bparticle{Van~den} \bsnm{Driessche}, \binits{P.}},
\bauthor{\bsnm{Watmough}, \binits{J.}}:
\bctitle{Further notes on the basic reproduction number}.
In: \bbtitle{Mathematical Epidemiology},
pp. \bfpage{159}--\blpage{178}.
\bpublisher{Springer},
\blocation{Berlin, Heidelberg}
(\byear{2008})
\end{bchapter}
\endbibitem

\bibitem{Riano}
\begin{botherref}
\oauthor{\bsnm{Ria{\~n}o}, \binits{G.}}:
Epidemic models with random infectious period.
medRxiv
(2020)
\end{botherref}
\endbibitem

\bibitem{busenberg1990analysis}
\begin{barticle}
\bauthor{\bsnm{Busenberg}, \binits{S.}},
\bauthor{\bparticle{van~den} \bsnm{Driessche}, \binits{P.}}:
\batitle{Analysis of a disease transmission model in a population with varying
  size}.
\bjtitle{Journal of mathematical biology}
\bvolume{28}(\bissue{3}),
\bfpage{257}--\blpage{270}
(\byear{1990})
\end{barticle}
\endbibitem

\bibitem{busenberg1993method}
\begin{barticle}
\bauthor{\bsnm{Busenberg}, \binits{S.}},
\bauthor{\bparticle{van~den} \bsnm{Driessche}, \binits{P.}}:
\batitle{A method for proving the non-existence of limit cycles}.
\bjtitle{Journal of mathematical analysis and applications}
\bvolume{172}(\bissue{2}),
\bfpage{463}--\blpage{479}
(\byear{1993})
\end{barticle}
\endbibitem

\bibitem{Der}
\begin{barticle}
\bauthor{\bsnm{Derrick}, \binits{W.}},
\bauthor{\bparticle{van~den} \bsnm{Driessche}, \binits{P.}}:
\batitle{A disease transmission model in a nonconstant population}.
\bjtitle{Journal of Mathematical Biology}
\bvolume{31}(\bissue{5}),
\bfpage{495}--\blpage{512}
(\byear{1993})
\end{barticle}
\endbibitem

\bibitem{Green97}
\begin{barticle}
\bauthor{\bsnm{Greenhalgh}, \binits{D.}}:
\batitle{Hopf bifurcation in epidemic models with a latent period and
  nonpermanent immunity}.
\bjtitle{Mathematical and Computer Modelling}
\bvolume{25}(\bissue{2}),
\bfpage{85}--\blpage{107}
(\byear{1997})
\end{barticle}
\endbibitem

\bibitem{LiGraef}
\begin{barticle}
\bauthor{\bsnm{Li}, \binits{M.Y.}},
\bauthor{\bsnm{Graef}, \binits{J.R.}},
\bauthor{\bsnm{Wang}, \binits{L.}},
\bauthor{\bsnm{Karsai}, \binits{J.}}:
\batitle{Global dynamics of a {SEIR} model with varying total population size}.
\bjtitle{Mathematical biosciences}
\bvolume{160}(\bissue{2}),
\bfpage{191}--\blpage{213}
(\byear{1999})
\end{barticle}
\endbibitem

\bibitem{SunHsieh}
\begin{barticle}
\bauthor{\bsnm{Sun}, \binits{C.}},
\bauthor{\bsnm{Hsieh}, \binits{Y.-H.}}:
\batitle{Global analysis of an {SEIR} model with varying population size and
  vaccination}.
\bjtitle{Applied Mathematical Modelling}
\bvolume{34}(\bissue{10}),
\bfpage{2685}--\blpage{2697}
(\byear{2010})
\end{barticle}
\endbibitem

\bibitem{Raz}
\begin{botherref}
\oauthor{\bsnm{Razvan}, \binits{M.}}:
Multiple equilibria for an {SIRS} epidemiological system.
arXiv preprint
(2001).
arXiv:0101051
\end{botherref}
\endbibitem

\bibitem{Arinovar}
\begin{barticle}
\bauthor{\bsnm{Yang}, \binits{W.}},
\bauthor{\bsnm{Sun}, \binits{C.}},
\bauthor{\bsnm{Arino}, \binits{J.}}:
\batitle{Global analysis for a general epidemiological model with vaccination
  and varying population}.
\bjtitle{Journal of Mathematical Analysis and Applications}
\bvolume{372}(\bissue{1}),
\bfpage{208}--\blpage{223}
(\byear{2010})
\end{barticle}
\endbibitem

\bibitem{Mart}
\begin{bbook}
\bauthor{\bsnm{Martcheva}, \binits{M.}}:
\bbtitle{An Introduction to Mathematical Epidemiology}
vol. \bseriesno{61}.
\bpublisher{Springer},
\blocation{New York, NY}
(\byear{2015})
\end{bbook}
\endbibitem

\bibitem{Thieme}
\begin{bbook}
\bauthor{\bsnm{Thieme}, \binits{H.R.}}:
\bbtitle{Mathematics in Population Biology}.
\bpublisher{Princeton University Press},
\blocation{Princeton, NJ}
(\byear{2018})
\end{bbook}
\endbibitem

\bibitem{Chavez}
\begin{bbook}
\bauthor{\bsnm{Brauer}, \binits{F.}},
\bauthor{\bsnm{Castillo-Chavez}, \binits{C.}},
\bauthor{\bsnm{Feng}, \binits{Z.}}:
\bbtitle{Mathematical Models in Epidemiology}.
\bpublisher{Springer},
\blocation{New York, NY}
(\byear{2019})
\end{bbook}
\endbibitem

\bibitem{KeMcK}
\begin{barticle}
\bauthor{\bsnm{Kermack}, \binits{W.O.}},
\bauthor{\bsnm{McKendrick}, \binits{A.G.}}:
\batitle{A contribution to the mathematical theory of epidemics}.
\bjtitle{Proc. R. Soc. Lond. Series A, Containing papers of a mathematical and
  physical character}
\bvolume{115}(\bissue{772}),
\bfpage{700}--\blpage{721}
(\byear{1927})
\end{barticle}
\endbibitem

\bibitem{earn2008light}
\begin{bchapter}
\bauthor{\bsnm{Earn}, \binits{D.J.}}:
\bctitle{A light introduction to modelling recurrent epidemics}.
In: \bbtitle{Mathematical Epidemiology},
pp. \bfpage{3}--\blpage{17}.
\bpublisher{Springer},
\blocation{Berlin, Heidelberg}
(\byear{2008})
\end{bchapter}
\endbibitem

\bibitem{Schaback}
\begin{barticle}
\bauthor{\bsnm{Schaback}, \binits{R.}}:
\batitle{On {COVID}-19 modelling}.
\bjtitle{Jahresbericht der Deutschen Mathematiker-Vereinigung}
\bvolume{122}(\bissue{3}),
\bfpage{167}--\blpage{205}
(\byear{2020})
\end{barticle}
\endbibitem

\bibitem{bacaer2020modele}
\begin{barticle}
\bauthor{\bsnm{Baca{\"e}r}, \binits{N.}}:
\batitle{Un mod{\`e}le math{\'e}matique des d{\'e}buts de l'{\'e}pid{\'e}mie de
  coronavirus en france}.
\bjtitle{Mathematical Modelling of Natural Phenomena}
\bvolume{15},
\bfpage{29}
(\byear{2020})
\end{barticle}
\endbibitem

\bibitem{Ketch}
\begin{botherref}
\oauthor{\bsnm{Ketcheson}, \binits{D.I.}}:
Optimal control of an {SIR} epidemic through finite-time non-pharmaceutical
  intervention.
arXiv preprint arXiv:2004.08848
(2020)
\end{botherref}
\endbibitem

\bibitem{Charp}
\begin{barticle}
\bauthor{\bsnm{Charpentier}, \binits{A.}},
\bauthor{\bsnm{Elie}, \binits{R.}},
\bauthor{\bsnm{Lauri{\`e}re}, \binits{M.}},
\bauthor{\bsnm{Tran}, \binits{V.C.}}:
\batitle{{COVID}-19 pandemic control: balancing detection policy and lockdown
  intervention under {ICU} sustainability}.
\bjtitle{Mathematical Modelling of Natural Phenomena}
\bvolume{15},
\bfpage{57}
(\byear{2020})
\end{barticle}
\endbibitem

\bibitem{Djidjou}
\begin{botherref}
\oauthor{\bsnm{Djidjou-Demasse}, \binits{R.}},
\oauthor{\bsnm{Michalakis}, \binits{Y.}},
\oauthor{\bsnm{Choisy}, \binits{M.}},
\oauthor{\bsnm{Sofonea}, \binits{M.T.}},
\oauthor{\bsnm{Alizon}, \binits{S.}}:
Optimal {COVID}-19 epidemic control until vaccine deployment.
medRxiv
(2020)
\end{botherref}
\endbibitem

\bibitem{Sofonea}
\begin{botherref}
\oauthor{\bsnm{Sofonea}, \binits{M.T.}},
\oauthor{\bsnm{Reyn{\'e}}, \binits{B.}},
\oauthor{\bsnm{Elie}, \binits{B.}},
\oauthor{\bsnm{Djidjou-Demasse}, \binits{R.}},
\oauthor{\bsnm{Selinger}, \binits{C.}},
\oauthor{\bsnm{Michalakis}, \binits{Y.}},
\oauthor{\bsnm{Alizon}, \binits{S.}}:
Epidemiological monitoring and control perspectives: application of a
  parsimonious modelling framework to the {COVID}-19 dynamics in {F}rance
(2020)
\end{botherref}
\endbibitem

\bibitem{alvarez2020simple}
\begin{botherref}
\oauthor{\bsnm{Alvarez}, \binits{F.E.}},
\oauthor{\bsnm{Argente}, \binits{D.}},
\oauthor{\bsnm{Lippi}, \binits{F.}}:
A simple planning problem for {COVID}-19 lockdown.
Technical report,
National Bureau of Economic Research
(2020)
\end{botherref}
\endbibitem

\bibitem{horstmeyer2020balancing}
\begin{botherref}
\oauthor{\bsnm{Horstmeyer}, \binits{L.}},
\oauthor{\bsnm{Kuehn}, \binits{C.}},
\oauthor{\bsnm{Thurner}, \binits{S.}}:
Balancing quarantine and self-distancing measures in adaptive epidemic
  networks.
arXiv preprint arXiv:2010.10516
(2020)
\end{botherref}
\endbibitem

\bibitem{di2020optimal}
\begin{botherref}
\oauthor{\bsnm{Di~Lauro}, \binits{F.}},
\oauthor{\bsnm{Kiss}, \binits{I.Z.}},
\oauthor{\bsnm{Miller}, \binits{J.}}:
Optimal timing of one-shot interventions for epidemic control.
medRxiv
(2020)
\end{botherref}
\endbibitem

\bibitem{Franco}
\begin{botherref}
\oauthor{\bsnm{Franco}, \binits{E.}}:
A feedback {SIR} ({fSIR}) model highlights advantages and limitations of
  infection-based social distancing.
arXiv preprint arXiv:2004.13216
(2020)
\end{botherref}
\endbibitem

\bibitem{baker2020reactive}
\begin{botherref}
\oauthor{\bsnm{Baker}, \binits{R.}}:
Reactive social distancing in a {SIR} model of epidemics such as {COVID}-19.
arXiv preprint arXiv:2003.08285
(2020)
\end{botherref}
\endbibitem

\bibitem{caulkins2020long}
\begin{barticle}
\bauthor{\bsnm{Caulkins}, \binits{J.}},
\bauthor{\bsnm{Grass}, \binits{D.}},
\bauthor{\bsnm{Feichtinger}, \binits{G.}},
\bauthor{\bsnm{Hartl}, \binits{R.}},
\bauthor{\bsnm{Kort}, \binits{P.M.}},
\bauthor{\bsnm{Prskawetz}, \binits{A.}},
\bauthor{\bsnm{Seidl}, \binits{A.}},
\bauthor{\bsnm{Wrzaczek}, \binits{S.}}:
\batitle{How long should the {COVID}-19 lockdown continue?}
\bjtitle{Plos one}
\bvolume{15}(\bissue{12}),
\bfpage{0243413}
(\byear{2020})
\end{barticle}
\endbibitem

\bibitem{caulkins2021optimal}
\begin{barticle}
\bauthor{\bsnm{Caulkins}, \binits{J.P.}},
\bauthor{\bsnm{Grass}, \binits{D.}},
\bauthor{\bsnm{Feichtinger}, \binits{G.}},
\bauthor{\bsnm{Hartl}, \binits{R.F.}},
\bauthor{\bsnm{Kort}, \binits{P.M.}},
\bauthor{\bsnm{Prskawetz}, \binits{A.}},
\bauthor{\bsnm{Seidl}, \binits{A.}},
\bauthor{\bsnm{Wrzaczek}, \binits{S.}}:
\batitle{The optimal lockdown intensity for {COVID}-19}.
\bjtitle{Journal of Mathematical Economics}
\bvolume{93},
\bfpage{102489}
(\byear{2021})
\end{barticle}
\endbibitem

\bibitem{bartlett2020deterministic}
\begin{bchapter}
\bauthor{\bsnm{Bartlett}, \binits{M.S.}}:
\bctitle{Deterministic and stochastic models for recurrent epidemics}.
In: \bbtitle{Contributions to Biology and Problems of Health},
pp. \bfpage{81}--\blpage{110}.
\bpublisher{University of California Press},
\blocation{Berkeley and Los Angeles}
(\byear{2020})
\end{bchapter}
\endbibitem

\bibitem{baxendale2011sustained}
\begin{barticle}
\bauthor{\bsnm{Baxendale}, \binits{P.H.}},
\bauthor{\bsnm{Greenwood}, \binits{P.E.}}:
\batitle{Sustained oscillations for density dependent markov processes}.
\bjtitle{Journal of mathematical biology}
\bvolume{63}(\bissue{3}),
\bfpage{433}--\blpage{457}
(\byear{2011})
\end{barticle}
\endbibitem

\bibitem{kurtz1978strong}
\begin{barticle}
\bauthor{\bsnm{Kurtz}, \binits{T.G.}}:
\batitle{Strong approximation theorems for density dependent markov chains}.
\bjtitle{Stochastic Processes and their Applications}
\bvolume{6}(\bissue{3}),
\bfpage{223}--\blpage{240}
(\byear{1978})
\end{barticle}
\endbibitem

\bibitem{Brit}
\begin{bbook}
\bauthor{\bsnm{Britton}, \binits{T.}},
\bauthor{\bsnm{Pardoux}, \binits{E.}},
\bauthor{\bsnm{Ball}, \binits{F.}},
\bauthor{\bsnm{Laredo}, \binits{C.}},
\bauthor{\bsnm{Sirl}, \binits{D.}},
\bauthor{\bsnm{Tran}, \binits{V.C.}}:
\bbtitle{Stochastic Epidemic Models with Inference}.
\bpublisher{Springer},
\blocation{Switzerland}
(\byear{2019})
\end{bbook}
\endbibitem

\bibitem{naasell2013influence}
\begin{barticle}
\bauthor{\bsnm{Naasell}, \binits{I.}}:
\batitle{The influence of immunity loss on persistence and recurrence of
  endemic infections}.
\bjtitle{Bulletin of mathematical biology}
\bvolume{75}(\bissue{11}),
\bfpage{2079}--\blpage{2092}
(\byear{2013})
\end{barticle}
\endbibitem

\bibitem{Kuehn}
\begin{barticle}
\bauthor{\bsnm{Jard{\'o}n-Kojakhmetov}, \binits{H.}},
\bauthor{\bsnm{Kuehn}, \binits{C.}},
\bauthor{\bsnm{Pugliese}, \binits{A.}},
\bauthor{\bsnm{Sensi}, \binits{M.}}:
\batitle{A geometric analysis of the {SIR}, {SIRS} and {SIRWS} epidemiological
  models}.
\bjtitle{Nonlinear Analysis: Real World Applications}
\bvolume{58},
\bfpage{103220}
(\byear{2021})
\end{barticle}
\endbibitem

\bibitem{de2019some}
\begin{botherref}
\oauthor{\bparticle{De~la} \bsnm{Sen}, \binits{M.}},
\oauthor{\bsnm{Nistal}, \binits{R.}},
\oauthor{\bsnm{Alonso-Quesada}, \binits{S.}},
\oauthor{\bsnm{Ibeas}, \binits{A.}}:
Some formal results on positivity, stability, and endemic steady-state
  attainability based on linear algebraic tools for a class of epidemic models
  with eventual incommensurate delays.
Discrete Dynamics in Nature and Society
\textbf{2019}
(2019)
\end{botherref}
\endbibitem

\bibitem{li2002qualitative}
\begin{barticle}
\bauthor{\bsnm{Li}, \binits{J.}},
\bauthor{\bsnm{Ma}, \binits{Z.}}:
\batitle{Qualitative analyses of {SIS} epidemic model with vaccination and
  varying total population size}.
\bjtitle{Mathematical and Computer Modelling}
\bvolume{35}(\bissue{11-12}),
\bfpage{1235}--\blpage{1243}
(\byear{2002})
\end{barticle}
\endbibitem

\bibitem{feng2000endemic}
\begin{barticle}
\bauthor{\bsnm{Feng}, \binits{Z.}},
\bauthor{\bsnm{Thieme}, \binits{H.R.}}:
\batitle{Endemic models with arbitrarily distributed periods of infection i:
  fundamental properties of the model}.
\bjtitle{SIAM Journal on Applied Mathematics}
\bvolume{61}(\bissue{3}),
\bfpage{803}--\blpage{833}
(\byear{2000})
\end{barticle}
\endbibitem

\bibitem{Feng}
\begin{barticle}
\bauthor{\bsnm{Feng}, \binits{Z.}}:
\batitle{Final and peak epidemic sizes for {SEIR} models with quarantine and
  isolation}.
\bjtitle{Mathematical Biosciences \& Engineering}
\bvolume{4}(\bissue{4}),
\bfpage{675}
(\byear{2007})
\end{barticle}
\endbibitem

\bibitem{AAK}
\begin{barticle}
\bauthor{\bsnm{Avram}, \binits{F.}},
\bauthor{\bsnm{Adenane}, \binits{R.}},
\bauthor{\bsnm{Ketcheson}, \binits{D.}}:
\batitle{A review of matrix {SIR} arino epidemic models}.
\bjtitle{MDPI}
\bvolume{11}(\bissue{4}),
\bfpage{89}--\blpage{104}
(\byear{2021})
\end{barticle}
\endbibitem

\bibitem{arino2020simple}
\begin{barticle}
\bauthor{\bsnm{Arino}, \binits{J.}},
\bauthor{\bsnm{Portet}, \binits{S.}}:
\batitle{A simple model for {COVID-19}}.
\bjtitle{Infectious Disease Modelling}
\bvolume{5},
\bfpage{309}--\blpage{315}
(\byear{2020})
\end{barticle}
\endbibitem

\bibitem{buonomo2015note}
\begin{barticle}
\bauthor{\bsnm{Buonomo}, \binits{B.}}:
\batitle{A note on the direction of the transcritical bifurcation in epidemic
  models}.
\bjtitle{Nonlinear Analysis: Modelling and Control}
\bvolume{20}(\bissue{1}),
\bfpage{38}--\blpage{55}
(\byear{2015})
\end{barticle}
\endbibitem

\bibitem{Hurtado}
\begin{barticle}
\bauthor{\bsnm{Hurtado}, \binits{P.J.}},
\bauthor{\bsnm{Kirosingh}, \binits{A.S.}}:
\batitle{Generalizations of the ‘linear chain trick': incorporating more
  flexible dwell time distributions into mean field {ODE} models}.
\bjtitle{Journal of mathematical biology}
\bvolume{79}(\bissue{5}),
\bfpage{1831}--\blpage{1883}
(\byear{2019})
\end{barticle}
\endbibitem

\bibitem{KamSal}
\begin{barticle}
\bauthor{\bsnm{Kamgang}, \binits{J.C.}},
\bauthor{\bsnm{Sallet}, \binits{G.}}:
\batitle{Computation of threshold conditions for epidemiological models and
  global stability of the disease-free equilibrium {(DFE)}}.
\bjtitle{Mathematical biosciences}
\bvolume{213}(\bissue{1}),
\bfpage{1}--\blpage{12}
(\byear{2008})
\end{barticle}
\endbibitem

\bibitem{bacaer2021mathematiques}
\begin{botherref}
\oauthor{\bsnm{Baca{\"e}r}, \binits{N.}}:
Math{\'e}matiques et {\'e}pid{\'e}mies.
Cassini
(2021)
\end{botherref}
\endbibitem

\bibitem{lotka1939analyse}
\begin{bbook}
\bauthor{\bsnm{Lotka}, \binits{A.J.}}:
\bbtitle{Analyse D{\'e}mographique Avec Application Particuli{\`e}re {\`A}
  L'esp{\`e}ce humaine}.
\bpublisher{Hermann},
\blocation{Paris}
(\byear{1939})
\end{bbook}
\endbibitem

\bibitem{dietz1993estimation}
\begin{barticle}
\bauthor{\bsnm{Dietz}, \binits{K.}}:
\batitle{The estimation of the basic reproduction number for infectious
  diseases}.
\bjtitle{Statistical methods in medical research}
\bvolume{2}(\bissue{1}),
\bfpage{23}--\blpage{41}
(\byear{1993})
\end{barticle}
\endbibitem

\bibitem{diekmann1990definition}
\begin{barticle}
\bauthor{\bsnm{Diekmann}, \binits{O.}},
\bauthor{\bsnm{Heesterbeek}, \binits{J.A.P.}},
\bauthor{\bsnm{Metz}, \binits{J.A.}}:
\batitle{On the definition and the computation of the basic reproduction ratio
  {R0} in models for infectious diseases in heterogeneous populations}.
\bjtitle{Journal of mathematical biology}
\bvolume{28}(\bissue{4}),
\bfpage{365}--\blpage{382}
(\byear{1990})
\end{barticle}
\endbibitem

\bibitem{Van}
\begin{barticle}
\bauthor{\bparticle{van~den} \bsnm{Driessche}, \binits{P.}},
\bauthor{\bsnm{Watmough}, \binits{J.}}:
\batitle{Reproduction numbers and sub-threshold endemic equilibria for
  compartmental models of disease transmission}.
\bjtitle{Mathematical biosciences}
\bvolume{180}(\bissue{1-2}),
\bfpage{29}--\blpage{48}
(\byear{2002})
\end{barticle}
\endbibitem

\bibitem{Diek}
\begin{barticle}
\bauthor{\bsnm{Diekmann}, \binits{O.}},
\bauthor{\bsnm{Heesterbeek}, \binits{J.}},
\bauthor{\bsnm{Roberts}, \binits{M.G.}}:
\batitle{The construction of next-generation matrices for compartmental
  epidemic models}.
\bjtitle{Journal of the Royal Society Interface}
\bvolume{7}(\bissue{47}),
\bfpage{873}--\blpage{885}
(\byear{2010})
\end{barticle}
\endbibitem

\bibitem{perasso2018introduction}
\begin{barticle}
\bauthor{\bsnm{Perasso}, \binits{A.}}:
\batitle{An introduction to the basic reproduction number in mathematical
  epidemiology}.
\bjtitle{ESAIM: Proceedings and Surveys}
\bvolume{62},
\bfpage{123}--\blpage{138}
(\byear{2018})
\end{barticle}
\endbibitem

\bibitem{Shuai}
\begin{barticle}
\bauthor{\bsnm{Shuai}, \binits{Z.}},
\bauthor{\bparticle{van~den} \bsnm{Driessche}, \binits{P.}}:
\batitle{Global stability of infectious disease models using {L}yapunov
  functions}.
\bjtitle{SIAM Journal on Applied Mathematics}
\bvolume{73}(\bissue{4}),
\bfpage{1513}--\blpage{1532}
(\byear{2013})
\end{barticle}
\endbibitem

\bibitem{fall2007epidemiological}
\begin{barticle}
\bauthor{\bsnm{Fall}, \binits{A.}},
\bauthor{\bsnm{Iggidr}, \binits{A.}},
\bauthor{\bsnm{Sallet}, \binits{G.}},
\bauthor{\bsnm{Tewa}, \binits{J.-J.}}:
\batitle{Epidemiological models and {L}yapunov functions}.
\bjtitle{Mathematical Modelling of Natural Phenomena}
\bvolume{2}(\bissue{1}),
\bfpage{62}--\blpage{83}
(\byear{2007})
\end{barticle}
\endbibitem

\bibitem{roberts2007pluses}
\begin{barticle}
\bauthor{\bsnm{Roberts}, \binits{M.}}:
\batitle{The pluses and minuses of 0}.
\bjtitle{Journal of the Royal Society Interface}
\bvolume{4}(\bissue{16}),
\bfpage{949}--\blpage{961}
(\byear{2007})
\end{barticle}
\endbibitem

\bibitem{li2011failure}
\begin{botherref}
\oauthor{\bsnm{Li}, \binits{J.}},
\oauthor{\bsnm{Blakeley}, \binits{D.}}:
The failure of {R0}.
Computational and Mathematical Methods in Medicine
\textbf{2011}
\end{botherref}
\endbibitem

\bibitem{guo2006global}
\begin{barticle}
\bauthor{\bsnm{Guo}, \binits{H.}},
\bauthor{\bsnm{Li}, \binits{M.Y.}}:
\batitle{Global dynamics of a staged progression model for infectious
  diseases}.
\bjtitle{Mathematical Biosciences \& Engineering}
\bvolume{3}(\bissue{3}),
\bfpage{513}
(\byear{2006})
\end{barticle}
\endbibitem

\bibitem{ma2006generality}
\begin{barticle}
\bauthor{\bsnm{Ma}, \binits{J.}},
\bauthor{\bsnm{Earn}, \binits{D.J.}}:
\batitle{Generality of the final size formula for an epidemic of a newly
  invading infectious disease}.
\bjtitle{Bulletin of mathematical biology}
\bvolume{68}(\bissue{3}),
\bfpage{679}--\blpage{702}
(\byear{2006})
\end{barticle}
\endbibitem

\bibitem{Andr}
\begin{barticle}
\bauthor{\bsnm{Andreasen}, \binits{V.}}:
\batitle{The final size of an epidemic and its relation to the basic
  reproduction number}.
\bjtitle{Bulletin of mathematical biology}
\bvolume{73}(\bissue{10}),
\bfpage{2305}--\blpage{2321}
(\byear{2011})
\end{barticle}
\endbibitem

\bibitem{liu1986influence}
\begin{barticle}
\bauthor{\bsnm{Liu}, \binits{W.-M.}},
\bauthor{\bsnm{Levin}, \binits{S.A.}},
\bauthor{\bsnm{Iwasa}, \binits{Y.}}:
\batitle{Influence of nonlinear incidence rates upon the behavior of {SIRS}
  epidemiological models}.
\bjtitle{Journal of mathematical biology}
\bvolume{23}(\bissue{2}),
\bfpage{187}--\blpage{204}
(\byear{1986})
\end{barticle}
\endbibitem

\bibitem{liu1987dynamical}
\begin{barticle}
\bauthor{\bsnm{Liu}, \binits{W.-M.-}},
\bauthor{\bsnm{Hethcote}, \binits{H.W.}},
\bauthor{\bsnm{Levin}, \binits{S.A.}}:
\batitle{Dynamical behavior of epidemiological models with nonlinear incidence
  rates}.
\bjtitle{Journal of mathematical biology}
\bvolume{25}(\bissue{4}),
\bfpage{359}--\blpage{380}
(\byear{1987})
\end{barticle}
\endbibitem

\bibitem{georgescu2007global}
\begin{barticle}
\bauthor{\bsnm{Georgescu}, \binits{P.}},
\bauthor{\bsnm{Hsieh}, \binits{Y.-H.}}:
\batitle{Global stability for a virus dynamics model with nonlinear incidence
  of infection and removal}.
\bjtitle{SIAM Journal on Applied Mathematics}
\bvolume{67}(\bissue{2}),
\bfpage{337}--\blpage{353}
(\byear{2007})
\end{barticle}
\endbibitem

\bibitem{tang2008coexistence}
\begin{barticle}
\bauthor{\bsnm{Tang}, \binits{Y.}},
\bauthor{\bsnm{Huang}, \binits{D.}},
\bauthor{\bsnm{Ruan}, \binits{S.}},
\bauthor{\bsnm{Zhang}, \binits{W.}}:
\batitle{Coexistence of limit cycles and homoclinic loops in a {SIRS} model
  with a nonlinear incidence rate}.
\bjtitle{SIAM Journal on Applied Mathematics}
\bvolume{69}(\bissue{2}),
\bfpage{621}--\blpage{639}
(\byear{2008})
\end{barticle}
\endbibitem

\bibitem{plemmons1977m}
\begin{barticle}
\bauthor{\bsnm{Plemmons}, \binits{R.J.}}:
\batitle{{M}-matrix characterizations. {I}—nonsingular {M}-matrices}.
\bjtitle{Linear Algebra and its Applications}
\bvolume{18}(\bissue{2}),
\bfpage{175}--\blpage{188}
(\byear{1977})
\end{barticle}
\endbibitem

\bibitem{hyman1999differential}
\begin{barticle}
\bauthor{\bsnm{Hyman}, \binits{J.M.}},
\bauthor{\bsnm{Li}, \binits{J.}},
\bauthor{\bsnm{Stanley}, \binits{E.A.}}:
\batitle{The differential infectivity and staged progression models for the
  transmission of hiv}.
\bjtitle{Mathematical biosciences}
\bvolume{155}(\bissue{2}),
\bfpage{77}--\blpage{109}
(\byear{1999})
\end{barticle}
\endbibitem

\bibitem{AABGH}
\begin{botherref}
\oauthor{\bsnm{Avram}, \binits{F.}},
\oauthor{\bsnm{Adenane}, \binits{R.}},
\oauthor{\bsnm{Bianchin}, \binits{G.}},
\oauthor{\bsnm{Goreac}, \binits{D.}},
\oauthor{\bsnm{Halanay}, \binits{A.}}:
Stability analysis of a seven parameter sir-type model including loss of
  immunity, and disease and vaccination fatalities.
Preprint
(2021)
\end{botherref}
\endbibitem

\bibitem{horn2012matrix}
\begin{bbook}
\bauthor{\bsnm{Horn}, \binits{R.A.}},
\bauthor{\bsnm{Johnson}, \binits{C.R.}}:
\bbtitle{Matrix Analysis}.
\bpublisher{Cambridge university press},
\blocation{New York, NY}
(\byear{2012})
\end{bbook}
\endbibitem

\bibitem{Sen}
\begin{barticle}
\bauthor{\bparticle{De~la} \bsnm{Sen}, \binits{M.}},
\bauthor{\bsnm{Ibeas}, \binits{A.}}:
\batitle{On an {SE(Is)(Ih)AR} epidemic model with combined vaccination and
  antiviral controls for {COVID-19} pandemic}.
\bjtitle{Advances in Difference Equations}
\bvolume{2021}(\bissue{1}),
\bfpage{1}--\blpage{30}
(\byear{2021})
\end{barticle}
\endbibitem

\end{thebibliography}



\end{document}